\documentclass[11pt]{amsart}
\usepackage{graphicx}
\usepackage{amsmath}
\usepackage{amscd}
\usepackage{amssymb}
\usepackage{latexsym}
\usepackage{epstopdf}
\newcommand{\R}{\mathbb{R}}
\newcommand{\Q}{\mathbb{Q}}
\newcommand{\C}{\mathbb{C}}
\newcommand{\Z}{\mathbb{Z}}
\newcommand{\Proj}{\mathbb{P}}

\newcommand{\set}[1]{\{ #1 \}}

\newcommand{\norm}[1]{|| #1 ||}

\newcommand{\del}{\partial}
\newcommand{\delbar}{\overline{\partial}}

\newcommand{\Bl}[1]{\widetilde{#1}}

\newcommand{\pii}[1] {\frac{ #1 }{2\pi i}}
\newcommand{\Tor}{\mathbf{P}}

\newcommand{\ellip}[1]{\frac{#1\theta(#1-z)\theta'(0)}{\theta(#1)\theta(-z)}}
\newcommand{\ellipar}[2]{\frac{#1\theta(#1-(-#2+1)z)\theta'(0)}{\theta(#1)\theta(-(-#2+1)z)}}

\newcommand{\ellipinv}[1]{\frac{\theta(#1)\theta(-z)}{#1\theta(#1-z)\theta'(0)}}

\newcommand{\jacc}[2]{\frac{\theta(#1-#2 z)\theta(-z)}{\theta(#1-z)\theta(-#2 z)}}

\newcommand{\orbnormTh}[2]{\frac{\theta(#1-(-#2+1)z)\frac{\theta'(0)}{2\pi i}}{\theta(#1)\theta(-(-#2+1)z)}}

\newcommand{\orbellipar}[3]{\frac{#3\theta(#1-(-#2+1)z)\theta'(0)}{\theta(#1)\theta(-(-#2+1)z)}}

\newcommand{\ellnormTh}[1]{\frac{\theta(#1-z)\frac{\theta'(0)}{2\pi 
i}}{\theta(#1)\theta(-z)}}

\newcommand{\ep}{\varepsilon}
\newcommand{\twopi}[1]{\frac{#1}{2\pi i}}

\newtheorem{thm*}{Theorem}
\newtheorem{thm}{Theorem}
\newtheorem{prop}{Proposition}
\newtheorem{lem}{Lemma}
\newtheorem{rmk}{Remark}

\newtheorem{cor*}{Corollary}
\newtheorem{cor}{Corollary}
\newtheorem{ack*}{Acknowledgements}

\title{Equivariant Elliptic Genera and Local McKay Correspondences}
\author{Robert Waelder}
\begin{document}
\baselineskip = 14 pt
\parskip = 2 pt

\begin{abstract}
In this paper we prove an equivariant version of the McKay correspondence for the elliptic genus on open varieties with a torus action. As a consequence, we will prove the equivariant DMVV formula for the Hilbert scheme of points on $\C^2$.
\end{abstract}
\email{rwaelder@math.ucla.edu}
\maketitle

\section{Introduction}

Let $X$ be a smooth variety, and $G$ a finite group acting holomorphically on $X$. Suppose that the quotient $X/G$ possesses a crepant resolution of singularities $V$. The McKay correspondence refers to the identification of topological invariants of $V$ with orbifold analogues of these invariants associated to the action of $G$ on $X$. The classical example is the case in which $G \subset SU(2)$ is a finite subgroup acting on $\C^2$. Then the quotient $\C^2/G$ possesses a unique crepant resolution $V$, and the Euler characteristic of $V$ coincides with the orbifold Euler characteristic $e_{orb}(\C^2,G) = \frac{1}{|G|}\sum_{gh=hg}e((\C^2)^{g,h})$. Over the years, the McKay correspondence has exhibited a remarkable versatility towards generalization. In \cite{B} Batyrev investigated a more general class of resolutions $\Bl{X/G}\rightarrow X/G$, and proved the McKay correspondence for the $E$-function in this situation. More recently, Borisov and Libgober have proven a similar result for the elliptic genus \cite{BL}.

In this paper we will prove an equivariant analogue of the McKay Correspondence for the elliptic genus. The advantage of working in the equivariant setting is that, by localization, we can make sense of the elliptic genus even for open varieties. This allows us to prove a host of new formulas. One consequence of the work in this paper is a beautiful formula for the generating function of the equivariant elliptic genus of the Hilbert scheme of points on $\C^2$ (with the standard torus action):

$$\sum_{n>0} p^n Ell((\C^2)^{[n]};y,q,t_1,t_2) = 
\prod_{m\geq 0,n>0,\ell,k}\frac{1}{(1-p^n q^m y^\ell t_1^{k_1} t_2^{k_2})^{c(nm,\ell,k)}}.$$
The terms $c(m,\ell,k)$ are the coefficients in the expansion of $Ell(\C^2;y,q,t_1,t_2)$ in $y,q,t_1$, and $t_2$. The above formula is an equivariant generalization of the DMVV formula:
$$\sum_{n>0} p^n Ell(S^{[n]};y,q) = 
\prod_{m\geq 0,n>0,\ell,k}\frac{1}{(1-p^n q^m y^\ell)^{c(nm,\ell,k)}}.$$
In the above formula, $S$ is a compact algebraic surface, and $Ell(S^{[n]};y,q)$ is the elliptic genus of the Hilbert scheme of $n$ points on $S$.
The non-equivariant DMVV formula was conjectured by string theorists Dijkgraaf, Moore, Verlinde and Verlinde \cite{DMVV}, and proven by Borisov and Libgober \cite{BL}. The equivariant version is a conjecture of Li, Liu, and Zhou \cite{LLJ}.

\subsection{Background on the Elliptic Genus}
For $X$ a smooth complex manifold, the elliptic genus of $X$ is defined as:

\begin{align}\label{Ell Def}
Ell(X) = \int_X\prod\frac{x_j\theta(\twopi{x_j}-z,\tau)}{\theta(\twopi{x_j},\tau)}
\end{align}

The product is taken over the Chern roots of the holomorphic tangent bundle to $X$. $\theta(\cdot,\tau)$ is the Jacobi theta function, and $z$ represents a formal parameter. Setting $y = e^{2\pi iz}$ and $q = e^{2\pi i\tau}$, the elliptic genus may also be interpreted as the index of the following differential operator:

\begin{align}\label{operator}
y^{-d/2}\delbar\otimes\bigotimes_{n=1}^{\infty}
\Lambda_{-yq^{n-1}}T^*X\otimes\Lambda_{-y^{-1}q^n}TX
\otimes S_{q^n}T^*X\otimes S_{q^n}TX
\end{align}

The modular properties of the Jacobi theta function endow the elliptic genus with a rich amount of structure. For example, if $X$ is Calabi-Yau, then the elliptic genus is a weak Jacobi form as a function of $(z,\tau)\in \C\times\mathbb{H}$. If $X$ is Calabi-Yau and possesses a nontrivial torus action, Liu has shown that the modular properties of the equivariant index of \ref{operator} actually imply its rigidity \cite{L}. In addition to these properties, the elliptic genus encodes a large number of classical algebraic and topological invariants of the space. For example, letting $q\to 0$ in the expression for the elliptic genus produces $y^{-d/2}$ times the Hirzebruch $\chi_{-y}$ genus, whereas letting $y\to 1$ produces the Euler characteristic of the space.

In \cite{CLW}, \cite{BLsing}, Chin-Lung Wang, Borisov, and Libgober investigated the following relative version of the elliptic genus for pairs $(X,D)$, where $D=\sum_i a_iD_i$ is a smooth divisor with normal crossings and coefficients $a_i\neq 1$:

$$Ell(X,D)=\int_X\prod\frac{x_j\theta(\twopi{x_j}-z)\twopi{\theta'(0)}}{\theta(\twopi{x_j})\theta(-z)}
\prod_i\frac{\theta(\twopi{c_1(D_i)}-(-a_i+1)z)\theta(z)}{\theta(\twopi{c_1(D_i)}-z)\theta((-a_i+1)z)}$$
The modular properties of the Jacobi theta function imply that the relative elliptic genus satisfies the following change of variables formula for blow-up morphisms: If $f:\Bl{X}\rightarrow X$ is the blow-up of $X$ along a smooth base with normal crossings with respect to the components of $D$, and $\Bl{D}$ is the divisor on $\Bl{X}$ satisfying: $K_{\Bl{X}}+\Bl{D} = f^*(K_X+D)$, then $Ell(\Bl{X},\Bl{D}) = Ell(X,D)$. For $Z$ a $\Q$-Gorenstein variety with log-terminal singularities, and $X\rightarrow Z$ a resolution of singularities with exceptional locus a normal crossing divisor $D$, Borisov and Libgober define the singular elliptic genus of $Z$ to be the relative elliptic genus of $(X,D)$. The change of variable formula, together with the Weak Factorization Theorem \cite{W} implies that this definition is well-defined. Moreover, when $Z$ possesses a crepant resolution $V$, the singular elliptic genus of $Z$ is easily seen to coincide with the elliptic genus of $V$. 

The singular elliptic genus (and its orbifold analogue) plays a crucial role in Borisov and Libgober's proof of the McKay correspondence for elliptic genera. Its utility stems from the fact that it behaves well with respect to a large class of birational modifications. The added flexibility obtained from studying the singular elliptic genus allowed Borisov and Libgober to reduce their proof to calculations involving toroidal embeddings. Their approach is similar in spirit to that of Batyrev in \cite{B}, who proved the McKay correspondence for the $E$-function by using the change of variables formula from motivic integration to reduce the case to calculations on toric varieties.

When $X$ has a nontrivial torus action, we may define the equivariant elliptic genus of $X$ to be the equivariant index of the operator defined in \ref{operator}. By the index theorem \cite{AS}, this is the same as the integral in equation \ref{Ell Def} obtained by replacing the Chern roots of $TX$ with their equivariant analogues. Similarly, we may define an equivariant version of the relative elliptic genus by replacing appearances of $c_1(D_i)$ with their equivariant extensions. 

The bulk of this paper is devoted to proving the change of variables formula for the equivariant orbifold elliptic genus (this case subsumes the non-orbifold case). Once we establish the change of variables formula in this situation, the remaining steps in the proof of the equivariant McKay correspondence for the elliptic genus follow closely the steps given in \cite{BL}.

\subsection{Outline of the Proof}
In a recent preprint \cite{RW}, I proved the equivariant change of variable formula for blow-ups along complete intersections $W = D_1\cap...\cap D_k\subset X$. The idea was to interpret the blow-down $\Bl{X}\rightarrow X$ as a toroidal morphism. The stratification defined by the divisors $D_i$ determined the toroidal structure of $X$, whereas the stratification defined by the proper transforms of these divisors, together with the exceptional divisor, determined the toroidal structure of $\Bl{X}$. The comparison of the relative elliptic genera of the base space and its blow-up was ultimately reduced to a computation involving the combinatorics of the polyhedral complexes associated to the two toroidal embeddings. This idea was inspired by Borisov and Libgober's use of polyhedral complexes in \cite{BL} to compute the push-forward of the orbifold elliptic class under the global quotient map.

Later it became apparant that the proof given in \cite{RW} could be adapted to the case in which $X$ was a ``normal cone space", i.e., a fiber product of spaces $\Proj(F\oplus 1)$, where $F\rightarrow W$ was a holomorphic vectorbundle. The idea was that the Chern roots of the tautological quotient bundle $Q_F\rightarrow X$ should play the role of the ``divisors" in a polyhedral complex associated to $X$. Similarly, if $f:\Bl{X}\rightarrow X$ was the blow-up of $X$ along $W$ with exceptional divisor $E$, then the Chern roots of $f^*Q_F\otimes\mathcal{O}(-E)$, and of $\mathcal{O}(E)$ should behave like the ``divisors" of a polyhedral complex associated to $\Bl{X}$. In this paper, we refer to such polyhedral complexes associated to the data of Chern roots as ``twisted polyhedral complexes." The case of a general blow-up may be reduced to cases of this nature by using an equivariant version of deformation to the normal cone.

The breakdown of the sections in this paper are as follows: For generic specializations of the parameters $(z,\tau)$, the integrand of the equivariant elliptic genus is a power series in the equivariant parameters with differential form coefficients. In section \ref{PSLocalization} we discuss convergence issues related to power series of this type and put their corresponding cohomology theory on solid ground. In section \ref{Definitions} we define our principle objects of study; namely the equivariant orbifold elliptic class and its relative version. In sections \ref{Toric Varieties} and \ref{Polyhedral Complexes} we review some facts about the equivariant cohomology of toric varieties and discuss how it relates to computing equivariant push-forwards of toroidal morphisms. In section \ref{Toroidal Pushforward} we use these results to compute the pushforward of the orbifold elliptic class under a toroidal morphism which is birational to a quotient by a finite group. This result is the equivariant analogue of Lemma $5.4$ in \cite{BL}. Sections \ref{Deformation Normal Cone} to \ref{Blow Up Formula} are devoted to the proof of the equivariant change of variable formula. In section \ref{Deformation Normal Cone}, we prove an equivariant analogue of deformation to the normal cone, tailored specifically to handle cohomological data like the orbifold elliptic class. As stated above, this will allow us to reduce the proof of the change of variable formula to the case when $X$ is a normal cone space. In section \ref{Normal Cone Space} we prove for completeness a number of technical lemmas regarding spaces of this form. In section \ref{Twisted Polyhedral Complex} we introduce the twisted polyhedral complex for normal cone spaces. In \ref{Blow Up Formula} we apply the techniques from the preceding sections to prove the equivariant change of variables formula. Finally, in \ref{Equiv McKay} we prove the equivariant McKay correspondence for elliptic genera, and the equivariant DMVV formula.

\begin{ack*}\rm
I wish to thank my advisor Professor Kefeng Liu for introducing me to elliptic genera and for his constant support, as well as Professor Anatoly Libgober for his feedback and help with technical aspects of his work.
\end{ack*}



\section{Equivariant Cohomology and Power Series}\label{PSLocalization}
\subsection{Preliminaries on Equivariant Cohomology}
We begin by reviewing some basic aspects of equivariant cohomology. 
For a thorough reference on the subject see \cite{AB}.

Let $M$ be a smooth manifold and $T$ a torus acting smoothly on $M$. 
Let $e_1,\ldots,e_\ell$ form a 
basis for the Lie algebra of $T$ which is dual to the linear forms 
$u_1,\ldots,u_\ell$. Every $X \in \mathfrak{t}$ defines a vectorfield 
$X$ on $M$ by the formula $X(p) = 
\frac{d}{dt}|_{t=0}\mathrm{exp}(tX)\cdot p$. Define $\Omega^*_T(M)$ to 
be the ring of differential forms on $M$ which are annihilated by 
$\mathcal{L}_X$ for every $X \in \mathfrak{t}$. If we let 
$d_{\mathfrak{t}} = d+\sum_{\alpha=1}^{\ell}u_\alpha i_{e_\alpha}$, 
then $d_{\mathfrak{t}}$ defines an operator on $\Omega^*_T(M)\otimes \C[u_1,\ldots,u_\ell]$ and 
satisfies $d_{\mathfrak{t}}^2 = 0$. The Cartan model for equivariant 
cohomology is defined to be:

$$H^*_T(M)_{\mathrm{Cartan}} = \frac{\ker d_{\mathfrak{t}}}
{\mathrm{im} d_{\mathfrak{t}}}.$$

The translation of concepts from cohomology to equivariant cohomology 
is more or less routine. For example, a $T$-map $f: M\rightarrow N$ 
induces a pullback $f^*:H^*_T(Y)\rightarrow H^*_T(X)$ as in ordinary cohomology. Similarly, for 
any $E \in K_T(X)$ we may define equivariant 
characteristic classes of $E$ which are equivariant extensions of the ordinary characteristic classes.

If $p$ is a single point with trivial $T$-action, the equivariant 
map $\pi : M \rightarrow p$ induces a map $\pi^*:H^*_T(p)\rightarrow 
H^*_T(M)$. Since $H^*_T(p) = \C[u_1,\ldots,u_\ell]$, the 
map $\pi^*$ makes $H^*_T(M)$ into a $\C[u_1,\ldots,u_\ell]$-module. 
Define $H^*_T(M)_{loc} = 
H^*_T(X)\otimes_{\C[u_1,\ldots,u_\ell]}\C(u_1,\ldots,u_\ell)$. A 
fundamental result of the subject is the localization theorem:

\begin{thm}
    Let $\set{P}$ denote the set of $T$-fixed components of $M$. Then
    $H^*_T(M)_{loc} \cong \bigoplus_P H^*(P)\otimes 
    \C(u_1,\ldots,u_\ell)$.
\end{thm}

If $P$ is a fixed component of $M$, the normal bundle to $P$ splits as 
a sum over the characters of the $T$-action on the fibers: $N_P = 
\bigoplus_{\lambda}V_{\lambda}$. Let $n^i_{\lambda}$ denote the formal Chern roots 
of $V_{\lambda}$. If we identify the equivariant parameters 
$u_1,\ldots,u_\ell$ with linear forms on the Lie algebra of $T$, then
the equivariant Euler class $e(P)$ of $N_P$ is equal to 
$\prod_{\lambda}\prod_{i}(n^i_{\lambda}+\lambda)$. Since 
none of the characters $\lambda$ are equal to zero, we see that $e(P)$ 
is always invertible. In light of this fact, we can describe the above 
isomorphism more explicitly. The map $H^*_T(X)_{loc} \rightarrow 
\bigoplus_P H^*(P)\otimes\C(u_1,\ldots,u_\ell)$ is given by $\omega 
\mapsto \bigoplus_P \frac{i_P^*\omega}{e(P)}$, where $i_P :P
\hookrightarrow X$ is the inclusion map.

If $f : M \rightarrow N$ is a proper map of $T$-spaces, we have the 
equivariant analogue of the cohomological push-forward $f_*: 
H^*_T(M)\rightarrow H^*_T(N)$. As in the non-equivariant setting, 
$f_*$ satisfies the projection formula $f_* (f^*(\omega)\wedge\eta) = 
\omega\wedge f_*\eta$. The new feature in equivariant cohomology is 
that we have an explicit expression for the restriction of 
$f_*\omega$ to a fixed component in $N$. This is given by the 
functorial localization formula \cite{MirrorI} \cite{MirrorII}:

\begin{thm}
    Let $f: M\rightarrow N$ be a proper map of $T$-spaces. Let $P$ be 
    a fixed component of $N$ and let $\set{F}$ be the collection 
    of fixed components in $M$ which $f$ maps into $P$. Let $\omega 
    \in H^*_T(M)$. Then:
    
    $$\sum_F f_*\frac{i_F^*\omega}{e(F)} = 
    \frac{i_P^*f_*\omega}{e(P)}.$$
\end{thm}

\subsection{Power Series in Equivariant Cohomology}
For simplicity, we assume that $M$ has a $T = S^1$ action. Let $X$ be the vectorfield on $M$ induced by the action of $T$. Let $C(M)$ denote the ring of formal power series in $u$ with coefficients in $\Omega^*(M)^T\otimes\C$. Then $d_X = d-ui_X : C(M) \rightarrow C(M)$. We define $H^*(C(M)) = \ker d_X/ \hbox{im } d_X$. 

Multiplication by $u$ gives $H^*(C(M))$ the structure of a $\C[u]$-module. Given any $\C[u]$-module $A$, we define $A_{loc} = A\otimes \C(u)$. Here, we interpret rational polynomials of the form $\frac{1}{1+uf(u)}$ as convergent power series by expanding around $u = 0$. 

Let $g$ be a $T$-invariant metric on $M$, and define $\theta = g(X,\cdot)$. Let $F \subset M$ be the $T$-fixed locus, and let $\rho$ be a bump function identically equal to one outside a tubular neighborhood of $F$, and equal to zero inside a smaller tubular neighborhood. (Note that we allow $F$ to have different components of varying dimensions). Then $\rho (d_X\theta)^{-1}$ is a well-defined element of $C(M)_{loc}$. Furthermore, it is easy to see that for any closed form $\omega \in C(M)_{loc}$, $\omega - d_X( \omega\cdot\theta \cdot \rho (d_X\theta)^{-1})$ has compact support in a tubular neighborhood of $F$. It follows that every closed form in $C(M)_{loc}$ may be represented by a form with compact support in a tubular neighborhood of $F$. Let $C(\nu_F)_{c}$ denote the ring of formal power series in $u$ whose coefficients are $T$-invariant forms with compact support in the normal bundle $\nu_F$. 

Let $\pi_* :C(\nu_F)_{c,loc} \rightarrow C(F)_{loc}$ be the map $\sum \omega_n u^n \mapsto \sum \pi_*\omega_n u^n$, where $\pi_*\omega_n$ denotes the integral of $\omega_n$ over the fiber of $\nu_F$. We first remark that $\pi_*i_X \omega = 0$ for any differential form $\omega$ with compact support in $\nu_F$. To see this, note that we can always express $\omega$ locally in the form $\omega = f(x,t)\pi^*\phi dt_1...dt_k$. Here $x$ represents the coordinates along $F$, $t$ represents the coordinates along the fiber, and $\pi^*\phi$ is the pullback of a form on $F$ via the projection $\pi : \nu_F \rightarrow F$. Thus $i_X\omega = f(x,t)\pi^*\phi\cdot i_X dt_1...dt_k$ which clearly integrates to zero along the fiber, since the degree of the form along the fiber is necessarily smaller than the fiber dimension.

It follows that $\pi_* d_X = d\pi_*$, and therefore it induces a map in cohomology $\pi_* : H^*(C(\nu_F)_c)_{loc} \rightarrow H^*(C(F))_{loc}$.

\begin{thm}
$\pi_*$ is injective.
\end{thm}

\begin{proof}
Let $\omega = \sum_n \omega_nu^n \in C(\nu_F)_{c,loc}$ be closed. Suppose $\pi_*\omega = d\eta$. By subtracting off $\pi^*d\eta\cdot \Phi$, where $\Phi$ is the equivariant Thom class of $\nu_F$, we may reduce to the case in which $\pi_*\omega = 0$. For any differential form $\alpha$, denote by $\alpha[i]$ the degree $i$ part. Then for some $0\leq k \leq \dim M$, $\omega = \omega[k]+\omega[k-1]+...+\omega[0]$, where $\omega[i] = \sum\omega_n[i]u^n$. Since $\omega$ is $d_X$-closed, $d\omega_n[k]=0$. Since $\pi_*\omega[k] = \sum\pi_*\omega_n[k]u^n = 0$, by the Thom isomorphism theorem, we can find compactly supported forms $y_n$ such that $dy_n = \omega_n[k]$. By averaging over $T$, we may assume that these forms are $T$-invariant. Then $\omega - d_X \sum y_n u^{n-1}$ has top degree $< k$ and is annihilated by $\pi_*$. The proof then follows by induction.
\end{proof}

Since $[\omega] = 0$ if and only if $\pi_*[\omega] = 0$, we see that $[\omega] = [\pi^*\pi_*\omega\cdot \Phi]$. Hence $[\omega]|_F = [\pi_*\omega\cdot e(\nu_F)]$, where $e(\nu_F)$ is the equivariant Euler class. Since the Euler class is invertible, the restriction map $[\omega] \mapsto [\omega]|_F$ must be invertible. This proves the localization theorem for $H^*(C(M))_{loc}$:

\begin{thm}
Let $F \subset M$ denote the fixed locus of $T$. Then the restriction map $H^*(C(M))_{loc} \rightarrow H^*(C(F))_{loc}$ is an isomorphism.
\end{thm}

We now introduce the ring of analytic forms: Let $C^{an}(M) \subset C(M)$ denote the ring of forms $\sum \omega_nu^n$ with the property that the partial sums $\sum_{n=0}^{N}\omega_ns^n$ converge in the $C^{\infty}$ sense to a form $\omega \in \Omega^*(M)^T$ for $\norm{s}$ sufficiently small. Let $B^{an}(M) = d_X(C(M))\cap C^{an}(M)$. We define 

$$H^*(C^{an}(M)) = \frac{\ker d_X: C^{an}(M)\rightarrow C^{an}(M)}{B^{an}(M)}.$$

If $F$ denotes the $T$-fixed locus of $M$, the above proof of the localization theorem extends word for word to the case of analytic forms. 

For $s \in \C^*$, following Witten \cite{Witten}, let $d_s = d-si_X : \Omega^*(M)^T \rightarrow \Omega^*(M)^T$. Define $H^*_s(M) = \frac{\ker d_s}{\hbox{im }d_s}$. Again, the above proof of the localization theorem adapts easily to the case of $H^*_s(M)$. Define an equivalence relation $\sim$ on $\prod_{s\in \C^*}H^*_s(M)$ as follows: For $\omega, \omega' \in \prod_{s\in \C^*}H^*_s(M)$, we say that $\omega \sim \omega'$ if $\omega_s = \omega'_s$ for $\norm{s}$ sufficiently small. By $\omega_s$, of course, we mean the $s$-component of $\omega$. We denote the group $\prod_{s\in \C^*}H^*_s(M)/\sim$ by $W$.

\begin{prop}
There is a natural evaluation map $ev: H^*(C^{an}(M))\rightarrow W$ given by $[\sum\omega_nu^n] \rightarrow [\sum\omega_ns^n]$.
\end{prop}

\begin{proof}
Let $\omega = \sum\omega_nu^n \in C^{an}(M)$ be a closed form, and let $d_X \sum\eta_nu^n \in C^{an}(M)$. We wish to show that for $\norm{s}$ sufficiently small, $[\sum\omega_ns^n] = [\sum\omega_ns^n+d_s\sum\eta_ns^n]$ in $H^*_s(M)$. First, for $\norm{s}$ sufficiently small, both sums $\sum\omega_ns^n$ and $\sum d_s\eta_ns^n$ converge. Therefore, it suffices to prove that $\sum d_s\eta_ns^n$ is exact for these values of $s$. Let $F \subset M$ be the $T$-fixed locus. Then $\sum d_s\eta_n|_Fs^n = \sum d\eta_n|_F s^n.$ By Hodge theory, $\sum d\eta_n|_F s^n$ converges to an exact form. Hence $[\sum d_s\eta_n s^n]|_F = 0$, and therefore $[\sum d_s \eta_n s^n] = 0$ by localization.
\end{proof}

Let $F(x_1,...,x_r)$ be function which is holomorphic in a neighborhood of $(0,...,0)$. Let $F = \sum_{i_1,...,i_r}a_{i_1,...,i_r}x_1^{i_1}\cdots x_r^{i_r}$ be the corresponding power series expansion. Let $[\omega_1],...,[\omega_r]$ be equivariant forms in $H^*_T(M)$ such that $u | \deg_0 \omega_i$. Note that this property is independent of the choice of representatives for $[\omega_i]$. Then $F([\omega_1],...,[\omega_r])$ is a well-defined element of $H^*(C^{an}(M))$. We can see this as follows:

Fix representatives $\omega_1,...,\omega_r \in C^{an}(M)$. Write $\omega_i = \tilde{\omega_i}+f_iu$, where $f_i \in C^{\infty}(M)\otimes \C[u]$ and $\deg \tilde{\omega_i} > 0$. Then 

$$F(\omega_1,...,\omega_r) = \sum_{i_1,...,i_r}a_{i_1,...,i_r}(\tilde{\omega_1}+f_1u)^{i_1}\cdots
(\tilde{\omega_r}+f_ru)^{i_r}$$

$$= \sum_{j_1,...,j_r}\tilde{\omega_1}^{j_1}\cdots\tilde{\omega_r}^{j_r}\sum_{i_1,...,i_r}a_{i_1,...,i_r}
{i_1 \choose j_1}\cdots {i_r\choose j_r} (f_1u)^{i_1-j_1}\cdots (f_ru)^{i_r-j_r}$$

$$= \sum_{j_1,...,j_r}\partial^{(j_1)}\cdots\partial^{(j_r)}F(f_1u,...,f_ru)\tilde{\omega_1}^{j_1}\cdots\tilde{\omega_r}^{j_r}.$$

Here $\partial^{(k)} = \frac{\partial^k}{k!}$. The sum over $j_1,...,j_r \geq 0$ is finite because the forms $\tilde{\omega_i}$ are nilpotent. It is therefore clear that the above expression is a well-defined $d_X$-closed form in $C^{an}(M)$. We next show that the corresponding class in $H^*(C^{an}(M))$ is independent of the choice of generators for $[\omega_i]$. Let $H(s_1,t_1,...,s_r,t_r) = F(s_1+t_1,...,s_r+t_r)$. Then $H = F(s_1,...,s_r)+t_1\tilde{H_1}+...+t_r\tilde{H_r}$. Thus, $F(\omega_1+d_X\eta_1,...,\omega_r+d_X\eta_r) = H(\omega_1,d_X\eta_1,...,\omega_r,d_X\eta_r) = F(\omega_1,...,\omega_r)+d_X(\eta_1\tilde{H_1}+...+\eta_r\tilde{H_r})$, so $F([\omega_1],...,[\omega_r])$ is well-defined independent of our choice of generators. Before proceeding further, we point out an important property possessed by the forms of this type. 

Let $\widehat{C}(M) \subset C^{an}(M)$ denote the subring $\set{\hbox{closed invariant forms } \otimes \C\{u\}}$ where $\C\{u\}$ denotes the ring of power series in $u$ which converge in sufficiently small neighborhoods of the origin. Let $\widehat{H}(M) \subset H^*(C^{an}(M))_{loc}$ denote the subspace of forms $[\omega]$ such that $\omega|_F \in \widehat{C}(F)_{loc}$ for some representative $\omega \in [\omega]$. Then if $[\omega_1],...,[\omega_r] \in H^*_T(M)$, and $u | \deg_0 \omega_i$, then $F([\omega_1],...,[\omega_r]) \in \widehat{H}(M)$. This is because $\deg_0 \omega_i|_F \in \C[u]$ and $d\tilde{\omega_i}|_F = 0$. Note also that functorial localization holds in this situation, and implies that, for $f: M \rightarrow N$, $f_*\widehat{H}(M) \subset \widehat{H}(N)$. 

\begin{lem}
$ev : \widehat{H}(M) \rightarrow W(M)$ is injective.
\end{lem}

\begin{proof}
Let $[\omega] \in \widehat{H}(M)$, so that $\omega|_F = \omega_1f_1+...+\omega_kf_k$, where $\omega_i$ are closed forms on $F$ and $f_i$ are convergent Laurent series in $u$ in some small disk about the origin. If $[\omega] \neq 0$, then by localization, $[\omega]|_F \neq 0$. Without loss of generality, we may assume that $\omega_1,...,\omega_k$ are linearly independent as cohomology classes. We can always choose an arbitrarily small $s$ so that $f_1(s),...,f_k(s)$ are not all zero. We then have that $ev[\omega]_s|_F = f_1(s)[\omega_1]+...+f_k(s)[\omega_k] \neq 0$. Hence $ev[\omega] \neq 0 $ in $W(M)$.
\end{proof}

All the above machinary was put together to make the following argument: Let $f: M \rightarrow N$, and let $F = F([\omega_1],...,[\omega_r])$ for $[\omega_i] \in H^*_T(M)$, $u | \deg_0 \omega_i$. Suppose $\theta_n(s)$ is a sequence of $d_s$-closed forms which converge in the $C^{\infty}$ sense to a representative of $(ev F)_s$ for $\norm{s}$ sufficiently small. If we factor $f$ as an inclusion followed by a projection, then $f_*$ makes sense on the form level and $f_*\theta_n(s)$ converge to representatives of $(f_* ev F)_s$. Now suppose that $[f_*\theta_n(s)] = [b_n(s)]$ and that $b_n(s) \rightarrow b(s)$. Then $f_*\theta_n(s) = b_n(s) + d_s\eta_n(s)$, and since $f_*\theta_n(s)$ and $b_n(s)$ converge, we have that $[b(s)] = [f_* ev F]_s$. Now if $s \mapsto [b(s)]$ corresponds to a form $ev G$ for some $G([\alpha_1],...,[\alpha_k]), [\alpha_i] \in H^*_T(N)$, then $f_* ev F = ev f_* F = ev G$, and therefore $f_* F = G$. Thus, we can compute push-forwards of summations by applying the push-forward term-by-term, so long as the corresponding summation converges in the sense discussed in this section. This observation will be used implicitly throughout this paper; for example, it will allow us in sections \ref{Polyhedral Complexes} and \ref{Twisted Polyhedral Complex} to reduce computations involving convergent power series to computations involving polynomials.

\begin{rmk}\rm
In what follows we will be interested in the case where $T$ is a compact torus of arbitrary dimension. It is easy to see how to generalize the above machinary to this situation. The difficulty lies more in the notation than in any other aspect.
\end{rmk}

\section{The Orbifold Elliptic Class}\label{Definitions}

Let $X$ be a smooth projective variety with a holomorphic $T\times G$ 
action. Let $D = \sum_I \alpha_i D_i$ be a smooth $G$-normal 
crossing divisor with $\alpha_i < 1$. The $G$-normal condition means that
if $g \in \mathrm{stab}_G(x)$ and $x\in D_i$, then $g\cdot D_i = D_i$.
Let $X^{g,h}_\gamma$ be a 
connected component of the fixed locus $X^{g,h}$ for some commuting 
pair $g,h \in G$. Then $N_{X^{g,h}_\gamma/X}$ splits into character 
sub-bundles $\bigoplus_\lambda N_\lambda$, where $g$ (resp. $h$) acts 
on $N_\lambda$ as multiplication by $e^{2\pi i\lambda(g)}$ (resp. 
$e^{2\pi i\lambda(h)}$). Let $I^{g,h}_\gamma \subset I$ index the set 
of divisors $D^X_i$ which contain $X^{g,h}_\gamma$. Since $D$ is 
$G$-normal, $g$ (resp. $h$) acts on 
$\mathcal{O}(D_i)|_{X^{g,h}_\gamma}$ as multiplication by $e^{2\pi 
i\lambda_i(g)}$ (resp. $e^{2\pi i\lambda_i(h)}$) for every $i\in 
I^{g,h}_\gamma$. For $i \not \in I^{g,h}_\gamma$, we define 
$\lambda_i(g) = \lambda_i(h) = 0$.

Following \cite{BL}, we define the orbifold elliptic class associated to the pair 
$(X^{g,h}_\gamma,D)$ as follows: 
$\mathcal{E}ll_{orb}(X^{g,h}_\gamma,D)=$

\begin{align*}
    &(i_{X^{g,h}_\gamma})_*\bigg\{
    \prod_{TX^{g,h}_\gamma}\ellip{\frac{x_i}{2\pi i}}
    \prod_{\lambda,N_\lambda}\ellnormTh{\frac{x_{\lambda,c}}{2\pi i}
    +\lambda(g)-\lambda(h)\tau}\bigg\}\\
    &\prod_I\jacc{\frac{D_i}{2\pi i}+\lambda_i(g)-\lambda_i(h)\tau}
    {(-\alpha_i+1)}
\end{align*}

Here $x_i$ denote the equivariant Chern roots of $TX^{g,h}_\gamma$, 
$x_{\lambda,c}$ denote the equivariant Chern roots of $N_\lambda$, and 
(abusing notation)
$D_i$ denote the equivariant first Chern classes of the 
corresponding divisors. 

Finally we define the orbifold elliptic class of $(X,D,G)$ as follows:
$$\mathcal{E}ll_{orb}(X,D,G) = \frac{1}{|G|}\sum_{gh=hg,\gamma}\mathcal{E}ll_{orb}(X^{g,h}_\gamma,D).$$
When $G = 1$, we will write $\mathcal{E}ll(X,D)$ instead of $\mathcal{E}ll_{orb}(X,D,G)$ and refer to this object as the equivariant elliptic class of the pair $(X,D)$.
We view all such classes as an elements inside the ring
$\widehat{H}(X)$.

Let $\set{F}$ denote the collection of fixed components of $X$. For each $F$ let $e(F)$ denote the equivariant Euler class of the normal bundle to $F$. We define the equivariant orbifold elliptic index 

$$Ell_{orb}(X,D,G)\equiv \big (\frac{2\pi i\theta(-z)}{\theta'(0)}\big )^{\dim X}\sum_F \int_F\frac{\mathcal{E}ll_{orb}(X,D,G)}{e(F)}.$$
It is a convergent power series in the equivariant parameters which depends implicitly upon the value of the complex parameter $z$ and on the lattice parameter $\tau$ used in the definition of the Jacobi theta function. We define the equivariant elliptic index $Ell(X,D)$ similarly. 

\section{Toric Varieties and Equivariant Cohomology}\label{Toric Varieties}
For a good reference on toric varieties, see \cite{F}.
Let $X$ be a smooth complete toric variety of dimension $n$. We denote the fan of $X$ 
by $\Sigma_X$, the lattice of $X$ by $N_X$, and the big torus by $T_X$. Let $Y$ be a smooth 
complete toric variety which satisfies the following properties:

$(1)$: $N_X \subset N_Y$ is a finite index sublattice.

$(2)$: $\Sigma_X$ is a refinement of $\Sigma_Y$ obtained by adding 
finitely-many one dimensional rays.

There is an obvious map of fans $\nu :\Sigma_X \rightarrow \Sigma_Y$ 
which induces a smooth map $\mu : X \rightarrow Y$. We call a map 
induced by such a morphism of fans a \it{toric morphism}\rm. It is 
easy to verify that $\mu: T_X \rightarrow T_Y$ is a covering map with 
covering group $N_Y/N_X$. Thus, we may regard $Y$ as a $T_X$-space. 
Our goal in this section is to obtain a convenient description of the 
equivariant pushforward $\mu_* : H^*_T(X)\rightarrow H^*_T(Y)$ in 
terms of the combinatorics of $\Sigma_X$ and $\Sigma_Y$. Here $T = 
T_X$. 

We first note that fixed points $F$ of $X$ are in $1-1$ 
correspondence with $n$-dimensional cones $C_F \subset \Sigma_X$. 
Furthermore, the infinitesimal weights of the $T$-action on $N_F$ 
correspond to linear forms in $\mathrm{Hom}(N_X,\Z)$ which are dual 
to the generators of $C_F$ in $N_X$. With this identification in mind, for any $\omega \in H_T^*(X)$, the the collection of polynomials $\set{\omega|_F}_{F\in \hbox{Fix}(X)}$ defines a piecewise polynomial function on the fan $\Sigma_X$. Define $\C[\Sigma_X]$ to be the ring of all piecewise polynomial functions on the fan of $X$. It is well-known that the map $H_T^*(X)\rightarrow \C[\Sigma_X]$ described here is an isomorphism:

\begin{thm}
    $H^*_T(X) \cong \C[\Sigma_X]$.
\end{thm}

Via this identification, we define 
$\nu_*:\C[\Sigma_X]\rightarrow \C[\Sigma_Y]$ to be the map which makes 
the following diagram commute:

$$\begin{CD}
     \C[\Sigma_X]  @>\nu_*>>   \C[\Sigma_Y]\\
    @|                              @|\\
     H^*_T(X)      @>\mu_*>>   H^*_T(Y)\\ 
\end{CD}$$

Here we understand $\C[\Sigma_Y]$ to be the ring of piecewise 
polynomial functions on $\Sigma_Y$ with respect to the lattice $N_X$.

We now describe $\nu_*$ more explicitly. First notice that for $f \in 
\C[\Sigma_X]$, $\nu_*f$ is given by viewing $f|_F$ as the zero degree 
part of an equivariant cohomology class $\omega \in H^*_T(X)$, 
pushing $\omega$ forward by $\mu_*$, and then forming the piecewise 
polynomial function defined by the zero degree part of $\mu_*\omega$. 
Thus, let $C \subset \Sigma_Y$ be an $n$-dimensional cone. Let 
$\nu^{-1}C$ be the fan $\Sigma_C \subset \Sigma_X$ which is the union 
of $n$-dimensional cones $C_i$. Let $x^{C_i}_1,\ldots,x^{C_i}_n$ be 
the linear forms dual to $C_i$ and $x^C_1,\ldots,x^C_n$ the linear 
forms in $\mathrm{Hom}(N_Y,\Z) \subset \mathrm{Hom}(N_X,\Z)$ dual to $C$. 
By functorial localization:

$$(\nu_*f)_C = \sum_{C_i \subset \Sigma_X}f_{C_i}
\frac{\prod_{j=1}^{n} x^C_j}{\prod_{j=1}^{n} x^{C_i}_j}.$$

Similarly, we define $\nu^* : \C[\Sigma_Y]\rightarrow \C[\Sigma_X]$ to 
be the map which makes the following diagram commute:

$$\begin{CD}
     \C[\Sigma_Y]  @>\nu^*>>   \C[\Sigma_X]\\
    @|                              @|\\
     H^*_T(Y)      @>\mu^*>>   H^*_T(X)\\ 
\end{CD}$$

\begin{prop}
    $\nu^*(f) = f\circ \nu$
\end{prop}

\begin{proof}
    Let $\omega \in H^*_T(Y)$ be the form such that $\omega|_P = 
    f|_P$ for every fixed point $P$. Let $F \in \mu^{-1}(P)$. Then
    
    $$\begin{CD}
     H^*_T(Y)      @>\mu^*>>     H^*_T(X)\\
    @VVV                              @VVV\\
     H^*_T(P)      @>\mu_F^*>>   H^*_T(F)\\ 
\end{CD}$$
commutes. Hence $(\mu^*\omega)|_F = \mu_F^*(\omega|_P) = 
\mu_F^*(f_P) = f_P$. Thus $\nu^*(f)$ is the piecewise polynomial 
function which is equal to $f_{C_P}$ on every cone $C_F \in 
\nu^{-1}C_P$. This is precisely the piecewise polynomial $f\circ \nu$.
\end{proof}

The map $\nu^* : \C[\Sigma_Y] \rightarrow \C[\Sigma_X]$ makes 
$\C[\Sigma_X]$ into a $\C[\Sigma_Y]$-module. As such, we observe:

\begin{prop}
    $\nu_*$ is a $\C[\Sigma_Y]$-module homomorphism.
\end{prop}

\begin{proof}
    In other words, we wish to prove the projection formula 
    $\nu_*(f\nu^*g) = \nu_*(f)\cdot g$. This follows from 
    identifying $\nu_*$ with $\mu_*$, $\nu^*$ with $\mu^*$ and 
    invoking the projection formula from equivariant cohomology.
\end{proof}

\section{Toroidal Embeddings and Toroidal Morphisms}\label{Polyhedral Complexes}

\subsection{Definitions}
Let $X$ be a compact variety and $D_X = \sum_{I_X} D^X_i$ a divisor on $X$ 
whose irreducible components are smooth normal crossing divisors.
For $I\subset I_X$, let $X_{I,j}$ denote the $j$th connected component 
of $\cap_I D^X_i$. Let $X^{o}_{I,j} = X_{I,j}-\cup_{I^c}D^X_i$. The 
collection of subvarieties $X^{o}_{I,j}$ form a stratification of $X$.
Associated to these data is a polyhedral complex with integral 
structure defined as follows:

Corresponding to $X_{I,j}$, define 
$N_{I,j} = \Z e_{i_1,j}+\ldots+\Z e_{i_k,j}$ 
to be the free group on the elements $e_{i_1,j},\ldots,e_{i_k,j}$. 
Here $i_1,\ldots i_k$ are the elements of $I$. 
Define $C_{I,j}$ to be the cone in the first orthant of this lattice. 
Whenever $I'\subset I$ and $X_{I,j}\subset X_{I',j'}$ we have natural 
inclusion maps $N_{I',j'}\hookrightarrow N_{I,j}$ and 
$C_{I',j'}\hookrightarrow C_{I,j}$. Define $\Sigma_X$ to be the 
polyhedral complex with integral structure obtained by gluing the 
cones $C_{I,j}$ together according to these inclusion maps.

Let 
$\C[\Sigma_X]$ denote the ring of piecewise polynomial functions on 
$\Sigma_X$. Fix $C \subset \Sigma_X$. 
Define $f^C$ to be 
the piecewise polynomial function which is equal to 
$\prod_{j=1}^{\dim C}x^{C}_j$ on every cone containing $C$, and 
equal to zero everywhere else. As in the toric geometry case, there is a natural 
correspondence between piecewise linear functions on $\Sigma_X$ and 
Cartier divisors whose irreducible components are components of 
$D_X$. We denote the piecewise linear function corresponding to $D$ 
by $f^D$.

\subsection{Toroidal Morphisms}
Our primary interest in this section is the study of toroidal morphisms.
This is a map $\mu: (X,D_X,\Sigma_X) \rightarrow (Y,D_Y,\Sigma_Y)$ 
which satisfies the following:

$(1)$: $\mu : X-D_X \rightarrow Y-D_Y$ is an unramified cover.

$(2)$: $\mu$ maps the closure of a stratum in $X$ to the closure of a 
stratum in $Y$.

$(3)$: Let $U_y$ be an analytic neighborhood of $y \in Y$ such that the 
components of $D_Y$ passing through $y$ correspond to coordinate 
hyperplanes. Then for $x \in \mu^{-1}(y)$, there exists an analytic 
neighborhood $U_x$ of $x$ such that the components of $D_X$ passing 
through $x$ correspond to coordinate hyperplanes of $U_x$. Moreover, 
the map $U_x \rightarrow U_y$ is given by monomial functions in the 
coordinates. 

Corresponding to $\mu$, we can define a map $\nu : \Sigma_X \rightarrow \Sigma_Y$
as follows: Let $C_{I,i} \subset \Sigma_X$ and let $e_1,\ldots,e_k 
\in N_{I,i}$ be the generators of $C_{I,i}$ which correspond to the 
divisors $D^X_1,\ldots, D^X_k$. We have that 
$\mu(X_{I,i}) = Y_{J,j}$. Let $v_1,\ldots,v_\ell \in N_{J,j}$ be 
the generators of $C_{J,j}$ which correspond to the divisors 
$D^Y_1,\ldots,D^Y_\ell$. For $1 \leq s \leq k$, $1 \leq t \leq 
\ell$, define $a_{st}$ to be the coefficient of $D^X_s$ of the divisor 
$\mu^*(D^Y_t)$. Then we define $\nu(e_s) = \sum a_{st}v_t$. Note that 
if $(X,\Sigma_X) \rightarrow (Y,\Sigma_Y)$ is a smooth toric morphism 
of toric varieties, then $\nu :\Sigma_X \rightarrow \Sigma_Y$ is the 
natural morphism of polyhedral complexes.

We have the following proposition relating $\nu$ to $\mu$:

\begin{prop}\label{Axioms}
    If $C = C_{J,j} \subset \Sigma_Y$, then $\nu^{-1}C$ is the union 
    of fans $\Sigma_\alpha \subset \Sigma_X$ with the following 
    properties:
    
    $(1)$: $\Sigma_\alpha$ is a refinement of $C$ obtained by adding 
    finitely-many $1$-dim rays.
    
    $(2)$: The lattice $N_\alpha$ of $\Sigma_\alpha$ is a finite index 
    sub-lattice of $N_C$.
    
    $(3)$: The fans $\Sigma_\alpha$ are in $1-1$ correspondence with 
    connected components $U_\alpha$ of $\mu^{-1}(N_{Y_{J,j}^{o}})$. 
    The map $U_\alpha \rightarrow N_{Y_{J,j}^{o}}$ is a fibration 
    given by the smooth toric morphism 
    $\Tor_{\Sigma_\alpha,N_\alpha}\rightarrow \Tor_{C,N_C}$ along the 
    fiber, and a $d_\alpha = d(\Sigma_\alpha)$-cover of $Y_{J,j}^{o}$ along the base.
\end{prop}

For a proof, see \cite{BL}. 

\subsection{Pushforward formula for Polyhedral 
Complexes}\label{Poly}
Motivated by the description of the push-forward $\nu_*$ 
for toric morphisms, define $\nu_* : \C[\Sigma_X] \rightarrow 
\C[\Sigma_Y]$ as follows. Let $C \subset \Sigma_Y$ be an 
$n$-dimensional cone with dual linear forms $x^C_1,\ldots,x^C_n$. 
Then for $f \in \C[\Sigma_X]$, we define:

$$(\nu_*f)_C = \sum_\alpha d_\alpha \sum_{C_i \in \Sigma_\alpha}
f_{C_i}\cdot \frac{\prod_{j=1}^{n}x^C_j}{\prod_{j=1}^{n}x^{C_i}_j}$$

The second sum is taken over the cones $C_i \subset \Sigma_\alpha$ 
with the same dimension as $C$.

Let $V$ be the toric variety $\coprod_\alpha d_\alpha\cdot 
\Tor_{\Sigma_\alpha,N_\alpha}$ with polyhedral fan $\Sigma_V$. We 
have a natural toric morphism $V \rightarrow \C^n$. We can compactify 
$V$ and $\C^n$ to obtain a smooth toric morphism $\overline{V}\rightarrow 
\Proj^n$. If we view $f$ as a piece-wise polynomial function on the 
fan of $\overline{V}$, then the above formula simply corresponds to 
$(\nu_*f)_C$ where $\nu :\Sigma_{\overline{V}} \rightarrow 
\Sigma_{\Proj^n}$. This identification allows us to apply the tools of 
the previous section toward the study of $\nu_*$. 

We first observe that $(\nu_*f)_C$ is indeed a polynomial function. 
This follows from the above identification of $\nu_*$ with the 
equivariant pushforward of a toric morphism. Furthermore, if we 
define $\nu^* :\C[\Sigma_Y]\rightarrow \C[\Sigma_X]$ by the formula 
$\nu^*(f) = f\circ \nu$ then the projection formula:

$$\nu_*(f\nu^*g)_C = \nu_*(f)_C\cdot g_C$$
follows from the projection formula in equivariant cohomology.

\begin{prop}
    $\nu_*(f)$ is a piece-wise polynomial function.
\end{prop}

\begin{proof}
We first show that $\nu_*(f^C)$ is piece-wise polynomial.

Fix $f = f^C$. Suppose $\nu(C) \subset C_0$ for some $C_0 \subset \Sigma_Y$ of 
dimension $k = \dim C$. Then $\nu_*(f)_{C_0} = d(\Sigma_{C_0})\prod_{j=1}^{k}x^{C_0}_j$. 
Suppose $C_1$ is a cone 
containing $C_0$. We wish to show $(\nu_*f)_{C_1}$ is an extension 
of $(\nu_*f)_{C_0}$. 

Consider the toric morphism 
$\sigma:\Tor_{\Sigma_{C_1},N(\Sigma_{C_1})}\rightarrow \C^{\dim C_1}$
induced by the map $\nu: 
\Sigma_{C_1}\rightarrow C_1$. Let $D_1,\ldots, D_k$ be the divisors 
in $\Tor_{\Sigma_{C_1},N(\Sigma_{C_1})}$ which correspond to the 
generators of $C$. Then the piece-wise polynomial 
function $f \in \C[\Sigma_{C_1}]$ represents the equivariant Thom 
class of $D_1 \cap \dots \cap D_k$. Since $\sigma(D_1\cap\dots\cap D_k)$ 
is the affine subspace of $\C^{\dim C_1}$ corresponding to $C_0$, we 
have that $\sigma_*(f)$ is the degree of $\sigma$ along $D_1\cap\dots\cap 
D_k$ times the polynomial function which represents the equivariant 
Thom class of this subspace. But this implies that: 

$$\nu_*(f)_{C_1} =d(\Sigma_{C_1})\frac{[N(\Sigma_{C_1}):N(C_1)]}
{[N(\Sigma_{C_0}):N(C_0)]}\prod_{j=1}^{k}x^{C_0}_j = 
d(\Sigma_{C_0})\prod_{j=1}^{k}x^{C_0}_j.$$

We need to explain the last equality. If $C_0$ corresponds to the 
strata $Y^{o}_{I,j}$ and $U\rightarrow N_{Y^{o}_{I,j}}$ is the fibration 
in Proposition \ref{Axioms} corresponding to the subdivision $\Sigma_{C_0}$, 
then $d(\Sigma_{C_0})[N(\Sigma_{C_0}):N(C_0)]$ and
$d(\Sigma_{C_1})[N(\Sigma_{C_1}):N(C_1)]$ both give the number of 
points in the pre-image of a generic point in $N_{Y^{o}_{I,j}}$. 

Next suppose that $C$ is mapped to a cone $C_0$ of strictly larger dimension. 
Consider the toric morphism 
$\Tor_{\Sigma_{C_0},N(\Sigma_{C_0})}\rightarrow \Tor_{C_0,N(C_0)}$ 
induced by the map $\nu: \Sigma_{C_0}\rightarrow C_0$. The polynomial 
function $f \in \C[\Sigma_{C_0}]$ represents the Thom class of an 
exceptional toric subvariety. Thus $\nu_*(f) = 0$, and it is easy to 
verify that $\nu_*(f) = 0$ on every cone containing $C_0$. Thus, 
$\nu_*$ maps the elements $f^C$ to piecewise polynomial 
functions. Since these functions generate $\C[\Sigma_X]$ as 
a $\C[\Sigma_Y]$-module, the proposition follows from the projection 
formula.
\end{proof}

In what follows we assume that $\mu: X \rightarrow Y$ is an 
equivariant map of projective $T$-spaces. Furthermore, we assume that 
the irreducible components of $D_X$ and $D_Y$ are invariant under 
the $T$-action. Define a map $\rho_X: \C[\Sigma_X]\rightarrow 
H^*_T(X)$ as follows: Fix a cone $C = C_{I,i}$ which corresponds to 
a connected component of the intersection locus of the divisors 
$D_1,\ldots,D_k$. Define $\rho_X[f^C\cdot (f^{D_1})^{a_1}\dots 
(f^{D_k})^{a_k}] = \Phi_{X_{I,i}}\wedge 
D_1^{a_1}\wedge \ldots \wedge D_k^{a_k}$. Here 
$\Phi_{X_{I,i}}$ denotes the equivariant Thom class of $X_{I,i} \subset X$ 
and, by abuse of notation, $D_j$ denote the equivariant Thom classes 
of the divisors $D_j$.

\begin{lem}
    $\rho_X$ is a ring homomorphism.
\end{lem}

\begin{proof}
    Fix cones $C_1=C_{I_1,i_1}$ and $C_2=C_{I_2,i_2}$. It suffices to 
    prove the theorem for the polynomials $f^{C_1}$ and 
    $f^{C_2}$. Let $I=I_1\cup I_2$. Let $C_{I,i}$ denote the cones 
    which correspond to components of the intersection 
    $X_{I_1,i_1}\cap X_{I_2,i_2}$. Clearly 
    $$f^{C_1}f^{C_2} = \sum_{I,i}f^{C_{I,i}}\prod_{I_1\cap I_2}f^{D_j}.$$
    
    Thus $\rho_X(f^{C_1}f^{C_2}) = 
    \sum_{I,i}\Phi_{X_{I,i}}\prod_{I_1\cap I_2}D_j$. However, by the 
    equivariant version of the excess intersection formula, this is 
    precisely the formula for $\rho_X(f^{C_1})\rho_X(f^{C_2})$.
\end{proof}

\begin{lem}
    $\rho_X\nu^* = \mu^*\rho_Y$.
\end{lem}

\begin{proof}
    It suffices to check this for polynomials $f^{C_{I,k}}$.
    If $D$ is a divisor on $Y$ whose irreducible components are 
    components of $D_Y$, then $\nu^*f^D$ is the piecewise linear 
    function corresponding to $\mu^*D$. It follows that 
    $\rho_X\nu^*f^D = \mu^*\rho_Y f^D$. Since all the maps are ring 
    homomorphisms, this implies that $\rho_X\nu^* \prod_{j\in I} 
    f^{D_j} = \mu^*\rho_Y\prod_{j\in I}f^{D_j}$. Let $\mu^*D_i = 
    \sum_j a_{ij}E_j$ as Cartier divisors. As in the lemma in the 
    Appendix, choose equivariant Thom forms $\Phi_{E_j}$ and 
    $\Phi_{D_i}$ with support in small tubular neighborhoods of their 
    respective divisors so that:
    
    $$\mu^*\Phi_{D_i} = \sum_j a_{ij}\Phi_{E_j} +d\psi_i$$
    as forms. Here $\psi_i$ are equivariant forms with compact support 
    in $\mu^{-1}N_{D_i}$. Let $\set{I,k}$ index the connected 
    components of $\cap_I D_i$. If we choose $N_{D_i}$ sufficiently small, 
    then
    
    $$\prod_I \Phi_{D_i} = \sum_{I,k}(\prod_I\Phi_{D_i})_{I,k}$$
    where $(\prod_I\Phi_{D_i})_{I,k}$ is the extension by zero of the 
    form $\prod_I\Phi_{D_i}|_{N_{I,k}}$. 
    
    Now $\prod_I f^{D_i} = \sum f^{C_{I,k}}$ and clearly 
    $(\prod_I\Phi_{D_i})_{I,k}$ is a representative of 
    $\rho_Y(f^{C_{I,k}})$. We have that 
    
    $$\mu^*(\prod_I\Phi_{D_i})_{I,k} = \big\{\prod_I (\sum_j 
    a_{ij}\Phi_{E_j} + d\psi_i) \big\}_{\mu^{-1}N_{I,k}}$$
    where the subscript $\mu^{-1}N_{I,k}$ means the extension by zero 
    of the form restricted to this open set. Since the $\psi_i$ forms 
    have compact support in $\mu^{-1}N_{D_i}$, this form is 
    cohomologous to 
    
    $$\big\{ \prod_I \sum_j a_{ij}\Phi_{E_j} 
    \big\}_{\mu^{-1}N_{I,k}}.$$
    
    But this is in turn a representative of $\rho_X\nu^* f^{C_{I,k}}$.
\end{proof}

\begin{lem}\label{PushForward}
    $\mu_* \rho_X = \rho_Y\nu_*$.
\end{lem}

\begin{proof}
    Since $\rho_X\nu^* = \mu^*\rho_Y$ and the polynomials $f^C$ 
    generate $\C[\Sigma_X]$ as a $\C[\Sigma_Y]$-module, by the 
    projection formula it suffices to check $\mu_* \rho_X f^C = 
    \rho_Y\nu_*f^C$.
    
    Case $1$: $C_{I,i}$ is mapped by $\nu$ to a cone $C_{J,j}$ of the same dimension.
    
    From the proof of Proposition $4$, $\nu_*f^{C_{I,i}} = 
    df^{C_{J,j}}$ where $d$ is the degree of $\mu : X_{I,i}\rightarrow 
    Y_{J,j}$. Thus, $\rho_Y\nu_* f^{C_{I,i}} = d\Phi_{Y_{J,j}} = 
    \mu_*\nu_* f^{C_{I,i}}$.
    
    Case $2$: $C_{I,i}$ is mapped by $\nu$ into a cone of strictly 
    larger dimension.
    
    As shown in Proposition $4$, $\nu_*f^{C_{I,i}} = 0$, so 
    $\rho_Y\nu_*f^{C_{I,i}} = 0 = \mu_*\Phi_{X_{I,i}} = \mu_*\rho_X 
    f^{C_{I,i}}$.
\end{proof}

\begin{rmk}\rm
    It is clear using the formalism of section \ref{PSLocalization} that the above lemmas relating $\mu$ to $\nu$ extend 
    without difficulty to the ring $\C[[\Sigma_X]]$ of piecewise 
    convergent power series. In this situation, $\rho_X$ is a map 
    from $\C[[\Sigma_X]]$ into $\widehat{H}(X)$.
\end{rmk}

\section{Elliptic Genera and Toroidal Morphisms}\label{Toroidal Pushforward}

Let $(\hat{X},\sum_{I_{\hat{X}}}\hat{D}_j)$ and $(X,\sum_{I_X}D_i)$ be smooth projective toroidal embeddings 
with $T$-actions which leave $\hat{D}_j$ and $D_i$ 
invariant. Let $G$ be a finite group which acts toroidally on $\hat{X}$ and 
commutes with the action of $G$. Suppose that $\mu : \hat{X}\rightarrow 
X$ 
is a $T$-equivariant toroidal morphism which is birational to a global 
quotient by $G$. If $\alpha_i$ are coefficients of the 
irreducible components $D_i$, define $\beta_j$ so that
$\mu^*(K_{X}+\sum\alpha_i D_i) = K_{\hat{X}} +\sum\beta_j \hat{D}_j$. Then:

\begin{thm}\label{Toroidal McKay}
    $$\mu_*\mathcal{E}ll_{orb}(\hat{X},\sum\beta_j \hat{D}_j,G)
    =\mathcal{E}ll(X,\sum\alpha_i D_i).$$
\end{thm}

\begin{proof}
    The proof of this theorem follows almost word for word the proof 
    of Theorem $5.1$ in \cite{BL}. The only difference is that here we 
    examine the equivariant push-forward of equivariant cohomology 
    classes, whereas \cite{BL} examine the push-forward of 
    (non-equivariant) classes in the Chow ring. We reproduce the proof 
    here for completeness.
    
    We refer frequently to the notation in the previous section:
    For commuting elements $g,h \in G$, let $\hat{X}^{g,h}_\gamma$ 
    denote the $\gamma$-th fixed component of $(g,h)$. Since the 
    action of $G$ on $\hat{X}$ is toroidal, $\hat{X}^{g,h}$ may be 
    identified with $X_{I^{g,h},i}$ for some indexing set 
    $I^{g,h}\subset I_{\hat{X}}$.
    
    Consider the following class in $\widehat{H}(\hat{X})$:
    
    $$E = \frac{1}{|G|}\sum_{gh=hg; {\hat{X}}^{g,h}_\gamma}
    \Phi^T_{\hat{X}^{g,h}_\gamma}\cdot\prod_{I_{\hat{X}}-I^{g,h}_\gamma}
    \orbellipar{\pii{\hat{D}_j}}{\beta_j}{\pii{\hat{D}_j}}\cdot$$
    $$\prod_{I^{g,h}_\gamma}
    \orbnormTh{\pii{\hat{D}_j}+\epsilon_j(g)-\epsilon_j(h)\tau}{\beta_j}
    e^{2\pi i(-\beta_j+1)\epsilon_j(h)z}.$$
    
    Here $\Phi^T_{\hat{X}^{g,h}_\gamma}$ is the equivariant Thom 
    class of the subvariety $\hat{X}^{g,h}_\gamma$, and 
    $\epsilon_j(g)$, etc., are defined as in the definition of the 
    orbifold elliptic genus.
    
    Our first goal is to prove
    
    \begin{align}\label{Poly McKay}
	\mu_*E = \prod_{I_{X}}\orbellipar{\pii{{D}_i}}{\alpha_i}{\pii{{D}_i}}
    \end{align}
    
    To prove the above equality, we express both sides as the image 
    under $\rho_{\hat{X}}$ and $\rho_X$ of piece-wise convergent 
    power series, and apply the push-forward formula from the 
    previous section. To that end, let $F$ be the piece-wise convergent power series 
    defined as follows: Let $C^{g,h}_\gamma$ be the cone which 
    corresponds to $X^{g,h}_\gamma$. For each cone $C=C_{J,j}$ 
    containing $C^{g,h}_\gamma$, let $F^{g,h}_\gamma|_C=$
    
    $$\frac{1}{|G|}\prod_{J}\orbellipar{\frac{x^C_j}{2\pi 
  i}+\epsilon_j(g)-\epsilon_j(h)\tau}{\beta_j}{\frac{x^C_j}{2\pi i}}e^{2\pi 
  i(-\beta_j+1)\epsilon_j(h)z}.$$
  
     For cones $C$ not containing $C^{g,h}_\gamma$, we define 
     $F^{g,h}_\gamma|_C = 0$.
     In the above expression, $x^C_j$ are the piece-wise linear functions dual to the 
     generators of $C$. If $\hat{D}_j$ are the divisors which correspond to 
     the generators of $C$, then $\epsilon_j(g)$, etc., are the
     infinitesimal weights attached to the divisors, and $\beta_j$ 
     are the coefficients of $\hat{D}_j$. It is easy to see that $F^{g,h}_\gamma$ is a 
     well-defined piece-wise convergent power series, and that 
     $\rho_{\hat{X}}(F^{g,h}_\gamma)$ is the $X^{g,h}_\gamma$-th 
     summand in the expression for $E$. We therefore define $F = 
     \sum_{gh=hg,\gamma} F^{g,h}_\gamma$, so that 
     $\rho_{\hat{X}}(F) = E$.
     
     Similarly, define $H \in \C[[\Sigma_X]]$ to be the piece-wise 
     convergent power series:
     
     $$H|_C = \prod_{i=1}^{\dim C}\ellipar{\frac{x^C_i}{2\pi 
i}}{\alpha_i}.$$

     Clearly $\rho_X(H)$ is equal to the right-hand side of \ref{Poly 
     McKay}. We are therefore reduced to proving that $\nu_* F = H$.
     
     Let $C \subset \Sigma_X$ be a cone, and let 
     $\Sigma_\alpha$ be the subdivisions of $C$ lying in 
     $\Sigma_{\hat{X}}$ which get mapped to $C$ under $\nu$. Let 
     $N_\alpha$ denote the lattices of $\Sigma_\alpha$, and $N_C$ the 
     lattice of $\Sigma_C$. Referring to the notation of section
     \ref{Poly}, the formula for $\nu_*$ tells us that:
     
     $$(\nu_*F)_C = \sum_\alpha d_\alpha\sum_{C_j\subset \Sigma_\alpha}
     F|_{C_j}\frac{\prod_{i=1}^{\dim C} x^{C}_i}{\prod_{i=1}^{\dim C}
     x^{C_j}_i}.$$
     
     For each $C_j \subset \Sigma_\alpha$ with the same dimension as 
     $C$, note that $F^{g,h}_\gamma|_{C_j} \neq 0$ if and only if 
     $C^{g,h}_\gamma \subset C_j$, i.e., if and only if $g$ and $h$ 
     are elements of the group $G_\alpha = N_C/N_\alpha$. We may 
     therefore write the push-forward of $F$ as:
     
     \begin{align*}
	 &\prod_{i=1}^{\dim C}x^C_i\sum_{\alpha}\frac{d_\alpha}{|G|}
	 \sum_{g,h \in G_\alpha}\\ 
	 &\sum_{C_j\subset \Sigma_\alpha}
	 \frac{\theta(\twopi{x^{C_j}_i}+\epsilon_i(g)-\epsilon_i(h)\tau-(-\beta_i+1)z)\twopi{\theta'(0)}}
	 {\theta(\twopi{x^{C_j}_i}+\epsilon_i(g)-\epsilon_i(h)\tau)\theta(-(-\beta_i+1)z)}e^{2\pi i(-\beta_i+1)z}
     \end{align*}
     
     By lemma $8.1$ in \cite{BL}, this is equal to $\sum_\alpha 
     \frac{d_\alpha |G_\alpha|}{|G|}H|_C$. But since $\sum_\alpha 
     d_\alpha|G_\alpha|$ describes the number of points in the 
     pre-image of a generic point in a tubular neighborhood of the 
     closed stratum corresponding to $C$, the coefficient in front of 
     $H|_C$ in the above expression is $1$. This completes the proof of \ref{Poly McKay}. To complete the proof, we apply the projection formula together with the following lemma relating the Chern classes of $\hat{X}$ to those of $X$:
     
     \begin{lem}
     $$\frac{c(T\hat{X})_T}{\mu^*c(TX)_T} = 
     \frac{\prod_{I_{\hat{X}}}(1+c_1(\hat{D}_j)_T)}
     {\prod_{I_X}(1+c_1(D_i)_T)}$$
     in $H^*_T(\hat{X})_{loc}$
     \end{lem}
     
     For details, see \cite{BL}. The above lemma may be proved using an argument analogous to the proof 
     of lemma $5$ in \cite{RW}.
     \end{proof}

\section{Deformation to the Normal Cone}\label{Deformation Normal Cone}

If $Q$ is a holomorphic function in a neighborhood of the origin, 
then $Q$ determines a map $\varphi_Q : K_T(\cdot) \rightarrow 
\widehat{H}(\cdot)$ by the rule $E \mapsto \prod_i Q(e_i)$, 
where $e_i$ represent the equivariant Chern roots of $E$. If, in 
addition, $Q(0) = 1$, then $\varphi_Q$ is multiplicative in the sense 
that $\varphi_Q(E_1\oplus E_2) = \varphi_Q(E_1)\varphi_Q(E_2)$. 

More generally, let $H$ be a finite abelian group with characters $\set{\lambda}$. 
Let $\set{f_\lambda}$ be an assignment of a holomorphic function in a 
neighborhood of the origin for each character. Such a 
collection $\set{f_\lambda}$ determines a map $\psi:K_T(\cdot)\otimes 
R(H)\rightarrow \widehat{H}(\cdot)$ by the rule $\psi:E = \bigoplus_\lambda 
E_\lambda \mapsto \prod_{\lambda} \prod f_\lambda(e_{\lambda,i})$. Here 
$e_{\lambda,i}$ denote the equivariant Chern roots of the 
$T$-vectorbundle $E_\lambda$. Let us fix a multiplicative map 
$\varphi = \varphi_Q$ and a (possibly not multiplicative) map 
$\psi = \psi_{\set{f_\lambda}}$.  

Let $X$ be a compact $T\times H$-variety and let $Z$
be a smooth $T\times H$-invariant subvariety. Let $V$ be a connected 
component of $X^H$. In the proof of the lemma below, the only difficult case to examine is when $Z\cap 
V \equiv W$ is a proper subset of $V$. We therefore assume this throughout. Let $\pi: \Bl{X}\rightarrow 
X$ be the blow-up of $X$ along $Z$ and let $\Bl{V}$ be the proper 
transform of $V$. Clearly $N_{\Bl{V}/\Bl{X}}$ has the same 
$H$-character decomposition as $N_{V/X}$. We may therefore make 
sense of $\psi(N_{\Bl{V}/\Bl{X}})$. The goal of this section is to 
prove the following crucial technical lemma relating 
$\varphi(T\Bl{V})\psi(N_{\Bl{V}/\Bl{X}})$ to 
$\varphi(TV)\psi(N_{V/X})$.

\begin{lem}\label{deformation}
    $$\pi_* \varphi(T\Bl{V})\psi(N_{\Bl{V}/\Bl{X}}) -
    \varphi(T{V})\psi(N_{V/X}) = (i_W)_*\Theta.$$
    Here $i_W$ is the inclusion map and $\Theta \in \widehat{H}(W)$ is 
    a universal class which depends only on the data of $W$, 
    $N_{W/V}$, and 
    $i_W^*N_{Z/X}$. 
\end{lem}




The argument given below is an adaptation of the argument in \cite{CLW} 
which was given for the non-equivariant case with $\psi = 1$.
 
\begin{proof}
    Let $\Pi: M_X\rightarrow X\times\Proj^1$ be the blow-up along 
    $Z\times\set{\infty}$. We give $M_X$ the obvious $T\times H$ 
    action. Let $M_V$ be the proper transform of $V\times \Proj^1$. 
    Clearly $M_V \rightarrow V\times \Proj^1$ is the blow-up along 
    $W\times\set{\infty}$. It is easy to see that $N_{M_V/M_X} \in 
    K_T(M_V)\otimes R(H)$ has the same $H$-character decomposition as 
    $N_{V/X}$. 
    
    Define $N_0$ to be the sub-bundle of $i_W^*N_{Z/X}$ on which $H$-acts 
    trivially. Let $N_1 = i_W^*N_{Z/X}/N_0$. Then $i_W^*TX$ decomposes 
    as $TW\oplus N_{W/Z}\oplus N_0\oplus N_1$. Clearly $TW\oplus 
    N_0 = i_W^*TV$, and therefore $N_{V/W}=N_0$ is a sub-bundle of 
    $i_W^*N_{Z/X}$ with quotient $N_1$. 
    
    Let $g : M_X\rightarrow \Proj^1$ be the composition $M_X\rightarrow X\times 
    \Proj^1\rightarrow \Proj^1$. Then $\hbox{div}(g) = X 
    -\Bl{X}-\Proj(N_{Z/X}\oplus 1)$. Furthermore, 
    $\hbox{div}(g|_{M_V}) = V - \Bl{V}-\Proj(N_{W/V}\oplus 1)$. Let 
    $i_V$, $i_{\Bl{V}}$, and $i_{\Proj(N_{W/V}\oplus 1)}$ denote the 
    respective inclusion maps of these divisors in $M_V$. 
    
    Let $p: 
    \Proj(N_{W/V}\oplus 1)\rightarrow W$ be the obvious fibration, 
    and let $S$ denote the tautological bundle over $\Proj(N_{W/V}\oplus 1)$.
    
    \noindent \bf{CLAIM}\rm: 
    \begin{align}
       i_V^* N_{M_V/M_X} =&\hbox{ } N_{V/X} \\
       i_{\Bl{V}}^* N_{M_V/M_X} =&\hbox{ } N_{\Bl{V}/\Bl{X}}\\
       i_{\Proj(N_{W/V}\oplus 1)}^*N_{M_V/M_X} =&\hbox{ }
       p^*N_{W/Z}\oplus p^*N_1\otimes S^*
    \end{align}
  
    $(1)$ is obvious. To prove $(2)$, note first that 
    $TM_X|_{\Bl{V}}$ decomposes as $T\Bl{V}\oplus 
    N_{\Bl{V}/\Bl{X}}\oplus N_{\Bl{X}/M_X}|_{\Bl{V}}$ and also as 
    $T\Bl{V}\oplus N_{\Bl{V}/M_V}\oplus N_{M_V/M_X}|_{\Bl{V}}$. Then simply notice that in both decompositions, $N_{\Bl{V}/\Bl{X}}$ and $i_{\Bl{V}}^* N_{M_V/M_X}$ are the nontrivial $H$-eigenspaces. 
    
    
    
    To prove $(3)$, note that $TM_X|_{\Proj(N_{W/V}\oplus 1)}$ 
    decomposes in the following two ways:
    \begin{align*}
	&TM_X|_{\Proj(N_{W/V}\oplus 1)} =\\
        &T\Proj(N_{W/V}\oplus 1)\oplus\mathcal{O}(-1)_{\Proj(N_{W/V}\oplus 1)}\oplus 
        i_{\Proj(N_{W/V}\oplus 1)}^*N_{M_V/M_X} =\\
        &T\Proj(N_{W/V}\oplus 1)\oplus 
	N_{\Proj(N_{W/V}\oplus 1)/\Proj(N_{Z/X}\oplus 1)}
	\oplus i_{\Proj(N_{W/V}\oplus 1)}^*\mathcal{O}(-1)_{\Proj(N_{Z/X}\oplus 
	1)}\\
    \end{align*}
    Observing $i_{\Proj(N_{W/V}\oplus 
    1)}^*\mathcal{O}(-1)_{\Proj(N_{Z/X}\oplus 1)} = 
    \mathcal{O}(-1)_{\Proj(N_{W/V}\oplus 1)}$ it follows that
    $i_{\Proj(N_{W/V}\oplus 1)}^*N_{M_V/M_X} = 
    N_{\Proj(N_{W/V}\oplus 1)/\Proj(N_{Z/X}\oplus 1)}$. It is easy to 
    verify that this bundle in turn is equal to $p^*N_{W/Z}\oplus p^*N_1\otimes S^*$.
    
    Since $\hbox{div}(g|_{M_V}) = V - \Bl{V}-\Proj(N_{W/V}\oplus 1)$,
    as equivariant classes $V = \Bl{V}+\Proj(N_{W/V}\oplus 1)$. Let 
    $u$ be the equivariant Thom class of $\Bl{V}$ (that is, the Thom 
    class of its normal bundle), and let $v$ be the equivariant Thom 
    class of $\Proj(N_{W/V}\oplus 1)$. Then $u+v$ is the equivariant 
    Thom class of $V$. Since $V$ is disjoint from $\Bl{V}$ and
    $\Proj(N_{W/V}\oplus 1)$, we have the relations $u(u+v) = 
    v(u+v)=0$. Note also that $uv$ is the equivariant Thom class of 
    $\Proj(N_{W/V})$, which is the exceptional divisor of 
    $\Bl{V}\rightarrow V$.
    
    Let $f$ be the holomorphic function in a neighborhood of the 
    origin which satisfies the relation $Q(z) = \frac{z}{f(z)}$. Then 
    by the above claim:
    
    \begin{align*}
	&\varphi(TM_V)\psi(N_{M_V/M_X})f(u+v) =\hbox{ } 
	(i_V)_*\varphi(TV)\psi(N_{V/X}) \\
	&\varphi(TM_V)\psi(N_{M_V/M_X})f(u) =\hbox{ } 
	(i_{\Bl{V}})_*\varphi(T\Bl{V})\psi(N_{\Bl{V}/\Bl{X}}) \\
	&\varphi(TM_V)f(v) =\hbox{ }(i_{\Proj(N_{W/V}\oplus 1)})_*
	\varphi(T\Proj(N_{W/V}\oplus 1))\psi(N_{\Proj(N_{W/V}\oplus 
	1)/\Proj(N_{Z/X}\oplus 1)})\\
    \end{align*}
    
    Note that since $u$ and $v$ are equivariant Thom classes, their 
    degree zero part is an element of $C^{\infty}(M_V)\otimes 
    \mathfrak{t}^*$ which vanishes at the origin of $\mathfrak{t}$. 
    Hence $f(u)$, etc., are well-defined elements of $\widehat{H}(M_V)$.
	
    Since $f(z) = z+\ldots$, we can define $g = f^{-1}$ in a possibly 
    smaller neighborhood of the origin. Consider the two-variable 
    holomorphic function $h(z_1,z_2) = f(z_1+z_2).$ Clearly 
    $h(z_1,z_2) = h(g(f(z_1)),g(f(z_2)))$. Define $F(y_1,y_2) = 
    h(g(y_1),g(y_2))$. Then $F$ is holomorphic in a neighborhood of 
    the origin and $F(f(z_1),f(z_2))= f(z_1+z_2)$. From \cite{Hirz} we have 
    the formula:
    
    $$F(y_1,y_2)g'(y_1)g'(y_2) = \sum_{(i,j)\neq 
    (0,0)}\varphi(H_{ij})y_1^i y_2^j.$$
    
    Here $\varphi(H_{ij})$ is the non-equivariant genus induced by 
    $f$ of the degree $(1,1)$ hypersurface $H_{ij}\subset 
    \Proj^i\times \Proj^j$. If we plug in $f(u)$ and $f(v)$ for $y_1$ 
    and $y_2$, we get $F(f(u),f(v))g'(f(u))g'(f(v)) = 
    f(u+v)g'(f(u))g'(f(v))$. By the relations $u(u+v)=v(u+v)=0$, this 
    last term is cohomologous to $f(u+v)$. It is instructive to go 
    over this last point in detail. Write $f(x+y) = \sum a_n(x+y)^n$ 
    and $g'(f(x))g'(f(y)) = \sum b_{ij}x^iy^j$. Note that $b_{00} = 
    1$. Let $h_{NIJ} = \sum_{n\leq N} a_n(x+y)^n\cdot \sum_{i\leq 
    I,j\leq J}b_{ij}x^{i}y^{j}$. Then $h_{NIJ}\rightarrow 
    f(x+y)g'(f(x))g'(f(y))$ in the $C^{\infty}$ topology. Abusing 
    notation, let $u$ and $v$ denote fixed representatives for their 
    respective cohomology classes. We have 
    that $h_{NIJ}(u,v) = \sum_{n\leq N}a_n(u+v)^n + d\eta_{NIJ}$. 
    Evaluating at a sufficiently small $s$ gives: $h_{NIJ}(u,v)(s) = 
    \sum_{n\leq N}a_n(u(s)+v(s))^n + d\eta_{NIJ}(s)$. Taking the 
    limit in the indices $N,I,J$, we get $h(u,v)(s) = f(u(s)+v(s)) + 
    d\eta(s)$. We therefore have that $ev(h) = ev(f(u+v))$, and 
    therefore $h = f(u+v)$ in $\widehat{H}(M_V)$.
    
    At the same time, $f(u+v) = \sum_{(i,j)\neq 
    (0,0)}\varphi(H_{ij})f(u)^if(v)^j.$ Thus
    
    \begin{align*}
	&\varphi(TM_V)\psi(N_{M_V/M_X})f(u+v) =\hbox{ } 
	\varphi(TM_V)\psi(N_{M_V/M_X})f(u)+\\
	&\varphi(TM_V)\psi(N_{M_V/M_X})f(v)
    +\varphi(TM_V)\psi(N_{M_V/M_X})
    \sum_{i+j\geq 2}\varphi(H_{ij})f(u)^{i}f(v)^{j} \\
    \end{align*}
    
    From the relations $u^2 = -uv$ and $v^2 = -uv$, we have that 
    $f(u)^{i+1} = f(u)f(-v)^i$. Therefore, 
    
    \begin{align*}
	\sum_{i+j\geq 2}\varphi(H_{ij})f(u)^{i}f(v)^{j} =&\hbox{ }\\ 
	\sum_{i+j\geq 2, i\geq 1}&\hbox{ }\varphi(H_{ij})f(u)f(-v)^{i-1}f(v)^j
    + \sum_{j\geq 2}\varphi(H_{0j})f(v)^j \\
    \end{align*}
    
    We therefore have $\sum_{i+j\geq 
    2}\varphi(H_{ij})f(u)^{i}f(v)^{j} = uv\frac{f(u)}{u}G(v)+ vJ(v)$ for some 
    universal convergent power series $G$ and $J$. 
    Let $\nu = 
    i_{\Proj(N_{W/V}\oplus 1)}^*v$. Clearly $\nu = c_1(S)$. Let 
    $w = \nu|_{\Proj(N_{W/V})}$. Finally, for ease of notation, write 
    $N =p^*N_{W/Z}\oplus p^*N_1\otimes S^*$.
    
    \begin{align*}
	&\varphi(TM_V)\psi(N_{M_V/M_X})f(u+v) =\\ 
	&\varphi(TM_V)\psi(N_{M_V/M_X})f(u)+\varphi(TM_V)\psi(N_{M_V/M_X})f(v)+\\ 
	&\varphi(TM_V)\psi(N_{M_V/M_X})uv\frac{f(u)}{u}G(v)+
	\varphi(TM_V)\psi(N_{M_V/M_X})vJ(v)\\
    \end{align*}
    It follows that
    \begin{align*}
	&(i_V)_*\varphi(TV)\psi(N_{V/X})=\\ 
	&(i_{\Bl{V}})_*\varphi(T\Bl{V})\psi(N_{\Bl{V}/\Bl{X}})+
	(i_{\Proj(N_{W/V}\oplus 1)})_*
	\varphi(T\Proj(N_{W/V}\oplus 1))\psi(N)+\\ 
	&(i_{\Proj(N_{W/V})})_*\varphi(T\Proj(N_{W/V})\oplus S|_{\Proj(N_{W/V})})
	\psi(N|_{\Proj(N_{W/V})})G(w)+\\
	&(i_{\Proj(N_{W/V}\oplus 1)})_*\varphi(T\Proj(N_{W/V}\oplus 1)\oplus S)\psi(N)J(\nu)\\
    \end{align*}
    
    Now apply the push-forward $\Pi_*$ to the above equation. Note 
    that $\Pi\circ i_V$ is the inclusion $v \mapsto (v,0)$ in 
    $V\times \Proj^1$. $\Pi\circ i_{\Bl{V}}$ is the composition of the 
    blow-down $\Bl{V}\rightarrow V$ with the inclusion of $V$ at infinity 
    in $V \times \Proj^1$. From the blow-up diagram:
    
    $$\begin{CD}
    \Proj(N_{W/V}\oplus 1) @> i_{\Proj(N_{W/V}\oplus 1)} >> M_V \\
    @V\hat{\Pi}VV               @VV\Pi V \\
    W\times\set{\infty} @>> i_W > V\times\Proj^1 \\
    \end{CD}$$
    we have that $\Pi\circ i_{\Proj(N_{W/V}\oplus 1)} = i_W\circ\hat{\Pi}$. Finally, 
    $\Pi\circ i_{\Proj(N_{W/V})}$ is clearly the composition of the blow-down map 
    $\hat{\pi}: \Proj(N_{W/V})\rightarrow W$ with the inclusion $i_W$. Therefore, 
    applying the pushforward $\Pi_*$ gives the equation:
    
    \begin{align*}
	&\varphi(TV)\psi(N_{V/X})=\\ 
	&\pi_*\varphi(T\Bl{V})\psi(N_{\Bl{V}/\Bl{X}})+
	(i_W)_*\Big\{
	\hat{\Pi}_*\varphi(T\Proj(N_{W/V}\oplus 1))\psi(N)+\\ 
	&\hat{\pi}_*\varphi(T\Proj(N_{W/V})\oplus S|_{\Proj(N_{W/V})})
	\psi(N|_{\Proj(N_{W/V})})G(w)+\\
	&\hat{\Pi}_*\varphi(T\Proj(N_{W/V}\oplus 1)\oplus S)\psi(N)J(\nu)\Big\}\\
    \end{align*}
    
    Since the term in the curly braces depends only on the data of 
    $i_W^*N_{Z/X}$, $N_{W/V}$, and $W$, this proves the lemma.
\end{proof}

\begin{rmk}\rm
Of course in the above proof, if $Q(x) = \frac{\twopi{x}\theta(\twopi{x}-z)\twopi{\theta'(0)}}{\theta(\twopi{x})\theta(-z)}$ and $f_\lambda = \frac{e^{2\pi i\lambda(h)\tau z}\theta(\twopi{x}+\lambda(g)-\lambda(h)\tau-z)\twopi{\theta'(0)}}{\theta(\twopi{x}+\lambda(g)-\lambda(h)\tau)\theta(-z)}$, then $\varphi(TV)\psi(N_{V/X})$ is the equivariant elliptic class associated to the pair $(g,h)$ and the $(g,h)$-fixed component $V$.
\end{rmk}
\section{The Normal Cone Space}\label{Normal Cone Space}

Let $W \subset X$ be a connected $T$-invariant subvariety of a 
projective $T$-space $X$. Suppose the normal bundle $N_{W/X}$ splits into 
a composition $L_1\oplus\ldots\oplus L_k$ of $T$-vectorbundles. Define 
$p: X^* \rightarrow W$ to be the fiber bundle with fiber 
$p^{-1}(w) = \Proj(L_1\oplus 1)_w\times \ldots \times \Proj(L_k\oplus 
1)_w$. It 
is easy to see that $X^*$ contains a copy of $W$ with the same normal 
bundle $N_{W/X}$. In our proof of the blow-up formula, we will ultimately 
reduce all computations on $X$ to computations on the more manageable 
space $X^*$. We therefore devote this section to gathering some 
important facts about the topology of $X^*$.

If we give the trivial vectorbundle $1$ the trivial action, then the 
action of $T$ on $W$ lifts naturally to $X^*$. Give $L_i$ a metric so 
that $T$ acts on $L_i$ by isometries, and give $\Proj(L_i\oplus 1)$ 
the induced metric. 

Define vectorbundles $Q_i \rightarrow X^*$ as follows: For $w \in W$ 
and $(\ell_1,\ldots,\ell_k)$ lines in $(L_1\oplus 
1)_w,\ldots,(L_k\oplus 1)_w$, define 
$(Q_i)_{(w,\ell_1,\ldots,\ell_k)} = \ell_i^{\perp} \subset (L_i\oplus 
1)_w$. These bundles inherit natural $T$-actions. Observe furthermore 
that $i_W^*(Q_i)$ is naturally isomorphic to $L_i$.

Define $V_i \subset X^*$ to be the subvariety 
$\set{(w,\ell_1,\ldots,\ell_k): \ell_i = [0:1]}$. For $i = 
1,\ldots,k$, $V_i$ are connected $T$-subvarieties, with connected 
intersection locus, and $W = 
V_1\cap\ldots\cap V_k$.

Finally, let $f: \Bl{X^*}\rightarrow X^*$ denote the blow-up of $X^*$ along 
$Z = V_1\cap\ldots\cap V_j$ with 
exceptional divisor $E$. We have the following intersection-theoretic 
result:
\label{top Chern}
\begin{thm}
    $c_{top}(f^*Q_1\oplus\ldots\oplus f^*Q_j \otimes \mathcal{O}(-E))_T = 
    0.$
\end{thm}

\begin{proof}
    We will show that 
    $\mathrm{Hom}(L,f^*Q_1\oplus\ldots\oplus f^*Q_j)$ 
    has an equivariant global nowhere zero section, where $L$ is a 
    line bundle with the same equivariant first Chern class as 
    $\mathcal{O}(E)$.
    
    We first give an explicit construction of the line bundle $L$. 
    Let $0\leq t_1,\ldots, t_j$. Define 
    $S_{t_1,\ldots,t_j} \subset \Bl{X^*}$ to be the subset
    $\set{(w,[v_1:1],\ldots,[v_j:1],[v_1:\ldots:v_j],\overline{\ell}): 
    \norm{v_i}=t_i}$. Here $\overline{\ell}$ represents a point in 
    $\Proj(L_{j+1}\oplus 1)_w\times\ldots\times\Proj(L_k\oplus 1)_w$. 
    We will refer to points in $\Bl{X^*}$ which are not contained in 
    any $S_{t_1,\ldots,t_j}$ as \it points at infinity\rm.   
    Let $\rho : [0,\infty) \rightarrow \R$ be a bump function equal to 
    zero in the region $[0,1/3)$ and equal to one in the region 
    $[2/3,\infty)$. 
    For $0 \leq t_i \leq 1$ and $v = 
    (w,[v_1:1],\ldots,[v_j:1],[v_1:\ldots:v_j],\overline{\ell})$ a point in 
    $S_{t_1,\ldots,t_j}$, let $L_v \subset (L_1\oplus 1)_w\oplus\ldots\oplus 
    (L_j\oplus 1)_w$ be the span of the vector 
    $\tilde{v} =  
    ((1-\rho(t_1))v_1,\rho(t_1),\ldots,(1-\rho(t_j))v_j,\rho(t_j)).$ 
    Outside this set, we define 
    $L_v$ to be the span of the vector $0\oplus 1\oplus\ldots\oplus 
    0\oplus 1$. This clearly defines a smooth line bundle on $\Bl{X^*}$ 
    with the same equivariant first Chern class as $\mathcal{O}(E)$. 
    
    We now prove that $\mathrm{Hom}(L,f^*Q_1\oplus\ldots\oplus 
    f^*Q_j)$ has a global equivariant nowhere zero section. For $v_i 
    \in L_i$, define: 
    
    $$h_i(v_i) = \Biggl \{
    \begin{matrix} 
	(-v_i,\norm{v_i}^2) & \norm{v_i}\leq 1 \cr
        (-\frac{v_i}{\norm{v_i}^2},1) & \norm{v_i}\geq 1 \cr 
    \end{matrix}
    $$
    
    For $v = (w,[v_1:1],\ldots,[v_j:1],[v_1:\ldots:v_j],\overline{\ell})
    \in S_{t_1,\ldots,t_j}$, we define $L_v \rightarrow 
    (Q_1\oplus\ldots\oplus Q_j)_{f(v)}$ by $\tilde{v} \mapsto 
    (h_1(v_1),\ldots,h_j(v_j))$. This section extends in a natural way 
    to the points at infinity, giving us a \it continuous \rm nowhere 
    zero equivariant section $s_0$ of $\mathrm{Hom}(L,f^*Q_1\oplus\ldots\oplus 
    f^*Q_j)$. The section is only continuous because the function:
    
     $$h(x) = \Biggl \{
    \begin{matrix} 
	x^2 & |x|\leq 1 \cr
        \frac{1}{x^2} & |x| \geq 1 \cr 
    \end{matrix}
    $$
    is not smooth. However, we may remedy this by approximating our 
    continuous section $s_0$ by a smooth section and then averaging over 
    the group $T$. Since $s_0$ was nowhere 
    zero and fixed by the $T$-averaging process, the new smooth 
    section will remain nowhere zero after averaging over $T$.
\end{proof}

\begin{rmk}\rm
    The intuition behind the preceding theorem is that if $Z = 
    D_1\cap\ldots\cap D_j$ is the complete intersection of a 
    collection of normal crossing divisors, then the proper 
    transforms $\Bl{D}_i$ of $D_i$ are disjoint when we blow up along 
    $Z$. We 
    therefore have that 
    $c_{top}(\mathcal{O}(\Bl{D}_1)
    \oplus\ldots\oplus\mathcal{O}(\Bl{D}_j))_T= 0$. While the above theorem is known,
    the proof given here will be useful later.
\end{rmk}


One easily observes that $i_{V_i}^*Q_i = N_{V_i/X^*}$. We might 
expect, therefore, that $c_{top}(Q_i)_T$ was the equivariant Thom 
class of $V_i$. This is in fact the case, as the following lemma 
proves:

\begin{lem}\label{top Chern is Thom class}
    $c_{top}(Q_i)_T = {i_{V_i}}_*1$
\end{lem}

\begin{proof}
    By the naturality properties of the equivariant Chern classes, it 
    is enough to prove the above lemma for $X^* = 
    \Proj(L\oplus 1)$, where $L$ is a $T$-vectorbundle over $W$, and
    $Q$ is the universal quotient bundle over $\Proj(L\oplus 1)$. 
    Here, $W$ itself plays the role of $V_i$ in the statement of the 
    lemma. As above, we endow the trivial bundle $1$ with the trivial $T$-action.
    
    Let us first prove that $c_{top}(Q_i) = {i_{W}}_*1$ in the 
    non-equivariant category. Let $p: \Proj(L\oplus 1)\rightarrow W$ 
    denote the obvious projection map. Let $r = \mathrm{rk}(Q)$.
    From the exact sequence:
    
    $$0\to S\to p^*(L\oplus 1)\to Q \to 0$$
    we have that $c_{r}(Q) = \big (\frac{p^*c(L)}{1-c_1(S^*)}\big 
    )_{r}$, where $(\cdot)_r$ denotes the degree $r$ part. Let 
    $p_{w}: \Proj(L\oplus 1)_w \rightarrow w$ denote the restriction 
    of $p$ to the fiber over $w$. Then $c_r(Q)|_{\Proj(L\oplus 1)_w} = 
    c_1(S^*)^r|_{\Proj(L\oplus 1)_w}$, which clearly integrates to 
    $1$ over the fiber. 
    
    We next observe that $Q$ has a global no-where zero section away 
    from $W$. We define this section as follows: 
    $s: (w,[v:z]) \mapsto (w,[v:z],\frac{\overline{z}v}{\norm{v}^2},-1)$
    We therefore have that $c_r(Q)$ is exact away from $W$. Hence, by 
    subtracting off an exact form, we may represent $c_r(Q)$ by a 
    form which has compact support in a tubular neighborhood of $W$, 
    and which integrates to one along every fiber of the tubular 
    neighborhood. It follows that $c_{r}(Q) = {i_W}_*1$ at least in 
    the non-equivariant sense.
    
    In general, $c_r(Q)_T$ is at least an equivariant extension of 
    the Thom class of $W$. Moreover, $c_r(Q)_T$ is equivariantly 
    exact outside $W$. We prove this as follows: Observe first that the non-zero 
    section we defined on the complement of $W$ was in fact 
    equivariant. This section therefore induces a splitting $Q = 
    Q'\oplus \C$ outside of $W$. Since we endowed $1$ with the 
    trivial $T$-action, one may easily verify that the $\C$ in the 
    above splitting inherits a trivial action as well. It follows 
    that $c_r(Q)_T$ is equivariantly exact outside of $W$.
    
    By subtracting off an exact form, we get that $c_r(Q)_T$ may be 
    represented by an equivariant form which has compact support in a tubular 
    neighborhood of $W$, and which integrates to one along every 
    fiber of the tubular neighborhood. It follows that $c_r(Q)_T$ is 
    the equivariant Thom class of $W$.
\end{proof}

We next prove a formula relating the equivariant 
Chern class of the blow-up of $X^*$ along $Z$ to that of $X^*$. 
Note first that $TX^*$ splits holomorphically and equivariantly into 
a direct sum of sub-bundles $F\oplus M$, with $i_Z^*F = TZ$ and 
$i_Z^*M = N_{Z/X^*}$. We may therefore apply the following lemma to 
compare the equivariant Chern classes of $X^*$ and $\mathrm{Bl}_Z X^*$.

\begin{lem}
    Let $Y$ be a complex $T$-space, and $Z \subset Y$ a $T$-invariant 
    complex submanifold. Suppose that $TY$ splits holomorphically and 
    equivariantly into a direct sum of sub-bundles $F\oplus M$ such 
    that $i_Z^*F = TZ$ and $i_Z^*M = N_{Z/Y}$. Let $f: 
    \Bl{Y}\rightarrow Y$ be the blow-up of $Y$ along $Z$ with 
    exceptional divisor $E$. Then:
    
    \begin{align}
	c(T\Bl{Y})_T &= c(f^*F)_Tc(f^*M\otimes \mathcal{O}(-E))_Tc(\mathcal{O}(E))_T
    \end{align}
    in the ring $H^*_T(\Bl{Y})_{loc}$.
\end{lem}

\begin{proof}
    By localization, it suffices to prove the equality at every fixed 
    component in $\Bl{Y}$. Let $\Bl{P} \subset \Bl{Y}$ be a fixed 
    component which is the proper transform of a fixed component $P 
    \subset Y$. If $P$ is disjoint from $Z$, then the equality of 
    $(4)$ at $\Bl{P}\cong P$ is trivial. Otherwise, $\Bl{P}$ is equal 
    to the blow-up of $P$ at $P\cap Z$. Note that $i_{P}^*F$ 
    decomposes as $F_0\oplus F_1$, where $T$ acts trivially on $F_0$ and 
    nontrivially on $F_1$. Similarly, $i_P^*M = M_0\oplus M_1$. 
    Clearly $TP = F_0\oplus M_0$ and $N_{P/Y} = F_1\oplus M_1$.
    
    Applying $i_{\Bl{P}}^*$ to $(4)$, we have:
    
    \begin{align*}
	i_{\Bl{P}}^*(\mathrm{LHS}) =\hbox{ } &c(T\Bl{P})c(N_{\Bl{P}/\Bl{Y}})_T\\
	i_{\Bl{P}}^*(\mathrm{RHS}) =\hbox{ } 
	&c(f^*F_0)c(f^*M_0\otimes\mathcal{O}(-E))c(i_{\Bl{P}}^*\mathcal{O}(E))\\
	\hbox{ }&c(f^*F_1)_Tc(f^*M_1\otimes\mathcal{O}(-E))_T
    \end{align*}	
     
    Since $\Bl{P}$ is the blow-up of $P$ along $P\cap Z$ and $i_{P\cap 
    Z}^*M_0 = N_{P\cap Z/P}$, the relation
    $c(T\Bl{P}) = 
    c(f^*F_0)c(f^*M_0\otimes\mathcal{O}(-E))c(i_{\Bl{P}}^*\mathcal{O}(E))$ 
    is well-known (see \cite{FultonIntersection}). It suffices therefore to prove that
    $c(N_{\Bl{P}/\Bl{Y}})_T = 
    i_{\Bl{P}}^*(c(f^*F_1)_Tc(f^*M_1\otimes\mathcal{O}(-E))_T)$. To this end we 
    prove the following claim:
    
    \bf{CLAIM}: $N_{\Bl{P}/\Bl{Y}} \cong i_{\Bl{P}}^*(f^*F_1\oplus 
    f^*M_1\otimes\mathcal{O}(-E))$ \it as $T$-vectorbundles.\rm
    
    To prove this, consider $f$ as a map $f:\Bl{P}\rightarrow P$. For 
    simplicity of notation, write $E$ for $E\cap \Bl{P}$. View 
    $N_{\Bl{P}/\Bl{Y}}$ as a sheaf, i.e., $N_{\Bl{P}/\Bl{Y}}(U) = 
    \Gamma(U,N_{\Bl{P}/\Bl{Y}})$. The derivative $Df: 
    N_{\Bl{P}/\Bl{Y}}\rightarrow f^*N_{P/Y} = f^*F_1\oplus f^*M_1$ is 
    a sheaf map, and maps onto the subsheaf $f^*F_1\oplus f^*M_1(-E)$. 
    Here $f^*M_1(-E)$ represents the subsheaf of $f^*M_1$ 
    corresponding to sections of $f^*M_1$ which vanish along $E$. 
    Let $s_0$ denote the global section of $\mathcal{O}(E)$ induced by 
    the defining equations of $E$. Then $\otimes 
    s_0^{-1}:f^*M_1(-E)\rightarrow f^*M_1\otimes\mathcal{O}(-E)$ is a 
    sheaf isomorphism. Define $\beta = (id\oplus \otimes 
    s_0^{-1})\circ Df$. By computing in local coordinates, one 
    verifies easily that $\beta: N_{\Bl{P}/\Bl{Y}}\rightarrow 
    f^*F_1\oplus f^*M_1\otimes\mathcal{O}(-E)$ is a sheaf isomorphism. 
    Furthermore, for $s$ a local section of $N_{\Bl{P}/\Bl{Y}}$, 
    $s(p) = 0$ if and only if $(\beta s)(p) = 0$. Therefore, $beta$ 
    induces an isomorphism of the corresponding vectorbundles. 
    Moreover, since $\beta$ is clearly equivariant, the vectorbundle 
    isomorphism is equivariant.
    
    This completes the proof if $P\cap E < P$. Next, suppose that $\Bl{P} \subset E$.
    For this case, it suffices to prove the following fact:
    
    \bf{CLAIM} $i_E^*(T\Bl{Y}\oplus \C) = i_E^*(f^*F\oplus 
    f^*M\otimes\mathcal{O}(-E)\oplus\mathcal{O}(E))$.\rm
    
    We prove this as follows: Let $\pi: \Proj(N_{Z/Y})\rightarrow Z$ 
    be the natural projection, and let $S$ denote the tautological 
    bundle over $\Proj(N_{Z/Y})$. Then $i_E^*T\Bl{Y} = 
    T\Proj(N_{Z/Y})\oplus S$. 
    
    Since $i_E^*f^* = \pi^*i_Z^*$, we have that $i_E^*f^*F = 
    \pi^*i_Z^*F = \pi^*TZ$ and $i_E^*f^*M = \pi^*N_{Z/Y}$. Thus, 
    
    \begin{align*}
	i_E^*(f^*F\oplus 
        f^*M\otimes\mathcal{O}(-E)\oplus\mathcal{O}(E)) =
	&\pi^*TZ\oplus \pi^*N_{Z/Y}\otimes S^*\oplus S
    \end{align*}
    
    From the exact sequence $0\to S\to \pi^*N_{Z/Y}\to Q\to 0$, we 
    have that $\pi^*N_{Z/Y}\otimes S^* = \C \oplus Q\otimes S^*$, 
    where $Q$ is the tautological quotient bundle. The claim then 
    follows from the observation that $T\Proj(N_{Z/Y}) = 
    \pi^*TZ\oplus Q\otimes S^*$. This completes the proof of the lemma.
\end{proof}

\begin{rmk}\rm
    Note that if $\Bl{Y}$ is equivariantly formal, the above proof 
    implies that $(4)$ holds in the unlocalized ring $H^*_T(\Bl{Y})$.
\end{rmk}

We may rewrite the left-hand side of $(4)$ as 
$c(f^*TY)_Tc(f^*M)_T^{-1}c(f^*M\otimes\mathcal{O}(-E))_Tc(\mathcal{O}(E))_T$. 
Viewed as an element of $\widehat{H}(\Bl{Y})$, it is easy to verify that 
this expression remains the same if we replace $M$ by any bundle $M'$ 
with $i_Z^*M' = N_{Z/Y}$. We therefore obtain the following corollary 
pertaining to $X^*$:

\begin{cor} \label{blow up Chern class}
    Let $f:\Bl{X^*}\rightarrow X^*$ be the blow-up of $X^*$ along $Z = 
    V_1\cap\ldots\cap V_j$ for $j\leq k$ with exceptional divisor $E$. 
    Then the following formula holds in $\widehat{H}(\Bl{X^*})$:
    
    \begin{align*}
	c(T\Bl{X^*})_T = \frac{c(f^*TX^*)_T}
	{\prod_{i=1}^jc(f^*Q_i)_T}\prod_{i=1}^j
	c(f^*Q_i\otimes\mathcal{O}(-E))_T c(\mathcal{O}(E))_T
    \end{align*}
\end{cor}

We end this section with a technical lemma which is the blow-up 
analogue of lemma \ref{top Chern is Thom class}. 

\begin{lem}\label{blow up Thom class}
    For $1\leq i\leq k$, let $\Bl{V_i}$ be the proper transform of $V_i$ under the above 
    blow-up: $f: \Bl{X^*}\rightarrow X^*$. Then 
    $c_{top}(f^*Q_i\otimes\mathcal{O}(-E))_T$ is the equivariant Thom 
    class of $\Bl{V_i}$ for $1\leq i \leq j$, and $c_{top}(f^*Q_i)_T$ 
    is the equivariant Thom class of $\Bl{V_i}$ for $j+1\leq i\leq 
    k$. Moreover $i_{\Bl{V_i}}^*c(f^*Q_i\otimes\mathcal{O}(-E))_T = 
    c(N_{\Bl{V_i}/\Bl{X^*}})_T$ for $1\leq i \leq j$ and 
    $i_{\Bl{V_i}}^*c(f^*Q_i)_T = c(N_{\Bl{V_i}/\Bl{X^*}})_T$ for 
    $j+1\leq i \leq k$.
    
\end{lem}

\begin{proof}
    Let $1\leq i\leq j$. Recall the equivariant continuous nowhere vanishing section $s_0$ 
    of $\mathrm{Hom}(L,f^*Q_1\oplus\ldots\oplus f^*Q_k)$ 
    constructed in the proof of theorem \ref{top Chern}. Let $\pi_i$ 
    denote the projection from $f^*Q_1\oplus\ldots\oplus f^*Q_k$ onto 
    the $i$-th factor. Then $\pi_i\circ s_0$ is an equivariant section 
    of $f^*Q_i\otimes L^*$. From the construction of $s_0$, 
    it is clear that $s_0$ vanishes precisely along $\Bl{V_i}$. Hence, 
    the equivariant top Chern class of $f^*Q_i\otimes L^*$ 
    is the equivariant Thom class of $\Bl{V_i}$. Since $L$ has the 
    same equivariant first Chern class as $\mathcal{O}(E)$, this 
    proves $c_{top}(f^*Q_i\otimes\mathcal{O}(-E))_T = {i_{\Bl{V_i}}}_*1$.
    
    To prove the second part of the lemma for $1\leq i\leq j$, apply 
    $i_{\Bl{V_i}}^*$ to both sides of the equation in corollary 
    \ref{blow up Chern class}. The LHS becomes $c(T\Bl{V_i})_T
    c(N_{\Bl{V_i}/\Bl{X^*}})_T$ while the RHS becomes:
    
    \begin{align*}
	&\frac{c(f^*N_{V_i/X^*}
    \otimes\mathcal{O}(-E\cap\Bl{V_i}))_Tc(f^*TV_i)_T}
    {\prod_{\ell\neq i}^jc(f^*N_{V_\ell/X^*})_T}\times\\
    &\prod_{\ell\neq i}^jc(f^*N_{V_\ell/X^*}
    \otimes\mathcal{O}(-E\cap\Bl{V_i}))_T 
    c(\mathcal{O}(E\cap\Bl{V_i}))_T\\
    \end{align*}
	
    Here we have used the fact that $i_{V_\ell}^*Q_{\ell} = 
    N_{V_\ell/X^*}$. By corollary \ref{blow up Chern class}, the factor multiplying
    $i_{\Bl{V_i}}^*c(f^*Q_i\otimes\mathcal{O}(-E))_T = c(f^*N_{V_i/X^*}
    \otimes\mathcal{O}(-E\cap\Bl{V_i}))_T$ in the above expression is equal to 
    $c(T\Bl{V_i})_T$. We therefore have that 
    $i_{\Bl{V_i}}^*c(f^*Q_i\otimes\mathcal{O}(-E))_T = 
    c(N_{\Bl{V_i}/\Bl{X^*}})_T$. 
    
    Next, let $j+1\leq i\leq k$. From the proof of lemma \ref{top 
    Chern is Thom class} we know that $Q_i$ has an equivariant 
    nowhere zero section in the complement of $V_i$. Pulling back 
    this section by $f$ gives an equivariant nowhere zero section in 
    the complement of $f^{-1}V_i = \Bl{V_i}$. It follows that 
    $c_{top}(f^*Q_i)_T$ must be localized in a neighborhood of 
    $\Bl{V_i}$. By the equivariant Thom isomorphism, 
    $c_{top}(f^*Q_i)_T$ must be a complex number multiple of the 
    equivariant Thom class of $\Bl{V_i}$. But since 
    $f_*c_{top}(f^*Q_i)_T$ is the Thom class of $V_i$, this complex 
    number multiple must be equal to one. 
    
    The proof of the second part of the lemma for $j+1\leq i\leq k$ 
    is analogous to the proof given above for $1\leq i\leq j$.
\end{proof}

\section{Twisted Polyhedral Complex}\label{Twisted Polyhedral Complex}
Throughout, we assume the following:

(a) $X$ is a smooth projective variety with a holomorphic $T\times H$ 
action, where $T$ is a torus and $H$ is a finite abelian group.

(b) $Z \subset X$ is a $T$-invariant smooth subvariety.

(c) $V \subset X^H$ is a connected component.

(d) $W = V\cap Z$ is connected. (see remark at the end of this section)

As noted above, $i_W^*N_{Z/X}$ splits as $N_0\oplus N_1$, 
where $N_0$ denotes the sub-bundle of $i_W^*N_{Z/X}$ on which $H$ acts 
trivially. Furthermore, $N_1$ decomposes as a sum of sub-bundles 
$\oplus N_{\lambda}$ corresponding to the characters of the $H$-action 
on $N_1$. Finally, $N_{W/Z}$ also splits into character sub-bundles 
$\oplus N_\ep$.

Let $E_1,\ldots,E_\ell$, $D_1,\ldots,D_k$ be smooth normal crossing 
divisors on $X$ intersecting $Z$ normally and non-trivially. We label 
these divisors so that $i_W^*\mathcal{L}_{D_i} \hookrightarrow N_0$ and
$i_W^*\mathcal{L}_{E_j} \hookrightarrow N_1$.  
Write $i_W^*\mathcal{L}_{D_i} = \Delta_i$ and 
$i_W^*\mathcal{L}_{E_j} = \xi_j$. Write $N_0 = F_0\oplus 
\bigoplus\Delta_i$ and $N_1 = \bigoplus F_{\lambda}\oplus \bigoplus \xi_j$.

Define a new space $X^*$ as follows: $p: X^*\rightarrow W$ is the 
fiber bundle with fiber $p^{-1}(w) =$
\begin{align*}
    &\prod\Proj(N_\ep\oplus 1)_w\times \Proj(F_0\oplus 1)_w\times \\
    &\prod\Proj(\Delta_i\oplus 1)_w\times\prod\Proj(F_{\lambda}\oplus 
1)_w\times\prod\Proj(\xi_j\oplus 1)_w
\end{align*}

Clearly $X^*$ contains a copy of $W$ with normal 
bundle $N_{W/X}$. Moreover, as in the previous section, each of the 
bundles $N_\ep$, $F_{0}$, $F_{\lambda}$, $\Delta_i$, and $\xi_j$ are the restrictions to $W$ 
of global equivariant bundles $Q_{N_\ep}$, $Q_{F_0}$, 
$Q_{F_\lambda}$, $Q_{\Delta_i}$ and $Q_{\xi_j}$. Recall also that the 
top Chern classes of these bundles are Poincar\'e dual to varieties 
$V_{N_\ep}, V_{F_0}, V_{F_\lambda}, V_{\Delta_i}$, and $V_{\xi_j}$. 
The analogue of $Z$ in $X^*$ is the intersection locus
of the varieties $V_{F_0}, V_{F_\lambda},V_{\Delta_i},$ and 
$V_{\xi_j}$. We will continue to call this intersection locus $Z$.

We associate a polyhedral complex with integral structure 
$(\Sigma,N_\Sigma)$ to $X^*$ as follows: Let 
$s_\ep = \mathrm{rk}(N_\ep)$, 
$r = \mathrm{rk}(F_0)$, and let $r_{\lambda} = 
\mathrm{rk}(F_{\lambda})$. Let 
$w^{\ep,a},x^{b},y^{\lambda,c},d^i$, and $e^j$ be an integral basis for 
$N_\Sigma =\Z^{\sum_\ep s_\ep+r+\sum_{\lambda}r_{\lambda}+k+\ell}$, for $a = 
1,\ldots,s_\ep$, $b 
=1,\ldots,r$, $c=1,\ldots,r_{\lambda}$, $i=1,\ldots,k$, and $j=1,\ldots,\ell$. 
Define $\Sigma$ to be the cone in the first orthant of the
vectorspace $N_\Sigma\otimes\R$. Let $w_{\ep,a},x_{b},y_{\lambda,c},d_i$, and $e_j$ be the linear 
forms on this vectorspace which are dual to above basis vectors. 
Then $\C[\Sigma] = \C[w_{\ep,a},x_{b},y_{\lambda,c},d_i,e_j]$. Let $G = 
\prod_\ep S_{s_\ep}\times S_{r}\times \prod_{\lambda}S_{r_\lambda}$, where $S_{n}$ denotes the 
symmetric group in $n$ letters. Then $G$ acts on $\C[\Sigma]$ by 
permuting the linear forms $w_{\ep,a},x_{b},y_{\lambda,c}$ in the obvious manner. 
Consider the following correspondence:

\begin{align*}
    w_{\ep,a}& \longleftrightarrow T\hbox{-Chern roots of }Q_{N_\ep}\\
    x_{b}& \longleftrightarrow T\hbox{-Chern roots of }Q_{F_0}\\
    y_{\lambda,c}& \longleftrightarrow T\hbox{-Chern roots of }Q_{F_\lambda}\\
    d_i &\longleftrightarrow c_1(Q_{\Delta_i})_T\\
    e_j &\longleftrightarrow c_1(Q_{\xi_j})_T\\
\end{align*}

Such a correspondence defines a natural map 
$\rho:\C[\Sigma]^G\rightarrow H^*_T(X^*)$.

Let $\phi: \Bl{X^*}\rightarrow X^*$ denote the blow-up of $X^*$ 
along $Z$ with exceptional divisor $E$. Define $\Bl{\Sigma}$ to be the polyhedral complex obtained 
from $\Sigma$ by adding the ray through the vector $\sum x_{b}+\sum 
y_{\lambda,c}+\sum d_i+\sum e_j$.
As before, let $w_{\ep,a},x_{b},y_{\lambda,c},d_i$, and $e_j$ denote the 
linear forms on $\C[\Bl{\Sigma}]$ which are dual to vectors 
$w^{\ep,a},x^{b},y^{\lambda,c},d^i$, and $e^j$. Let $t$ be the linear 
form dual to $\sum x_{b}+\sum y_{\lambda,c}+\sum d_i+\sum e_j$. Then:

$$\C[\Bl{\Sigma}] \cong 
\C[t,w_{\ep,a},x_{b},y_{\lambda,c},d_i,e_j]/\prod_{b,c,\lambda,i,j}x_{b}y_{\lambda,c}d_i e_j.$$

$G$ acts on $\C[\Bl{\Sigma}]$ in the obvious manner. Consider the correspondence:

\begin{align*}
    w_{\ep,a}& \longleftrightarrow T\hbox{-Chern roots of }\phi^*Q_{N_\ep}\\
    x_{b}& \longleftrightarrow T\hbox{-Chern roots of }\phi^*Q_{F_0}\otimes\mathcal{O}(-E)\\
    y_{\lambda,c}& \longleftrightarrow T\hbox{-Chern roots of 
    }\phi^*Q_{F_\lambda}\otimes \mathcal{O}(-E)\\
    d_i &\longleftrightarrow c_1(\phi^*Q_{\Delta_i}\otimes\mathcal{O}(-E))_T\\
    e_j &\longleftrightarrow c_1(\phi^*Q_{\xi_j}\otimes\mathcal{O}(-E))_T\\
\end{align*}

By theorem \ref{top Chern}, this correspondence 
induces a well-defined homomorphism $\rho: \C[\Bl{\Sigma}]^G\rightarrow 
H^*_T(X^*)$. 

It is easy to see that $\nu^* : \C[\Sigma]^{G}\rightarrow 
\C[\Bl{\Sigma}]^{G}$ and similarly $\nu_* 
:\C[\Bl{\Sigma}]^{G} \rightarrow \C[\Sigma]^{G}$. We 
have the following important lemmas:

\begin{lem}\label{pushforward commutes}
    \begin{align}
	\rho\nu^* =& \phi^*\rho \label{a} \\
	\phi_*\rho =& \rho\nu_* \label{b}
    \end{align}
\end{lem}

\begin{proof}
    For notational convenience, let $\bar{w}_{\ep,a}, 
    \bar{x}_b,\ldots$ denote the $T$-Chern roots corresponding to 
    $w_{\ep,a}, x_b,\ldots$. With this notation, 
    $\rho(f)(w_{\ep,a},x_b,\ldots)$ is equal to 
    $f(\bar{w}_{\ep,a},\bar{x}_b,\ldots)$. We first prove \ref{a}: We 
    have: 
    \begin{align*}
	\phi^*\rho(f) =&\hbox{ } 
    f(\phi^*\bar{w}_{\ep,a},\phi^*\bar{x}_b,\ldots)\\
    =&\hbox{ }f(\bar{w}_{\ep,a},\bar{x}_b+\bar{t},\ldots) = 
    f(\rho(w_{\ep,a}),\rho(x_b+t),\ldots)\\
    =&\hbox{ }f(\rho\nu^*w_{\ep,a},\rho\nu^*x_b,\ldots) = 
    \rho\nu^*f\\
    \end{align*}
    
    \ref{a} allows us to reduce the proof of \ref{b} to the case where $f = 
    t^n$. Let $N = \hbox{ codim }Z$. For $n < N-1$, $\nu_*t^n = 0 = 
    \phi_*\rho(t)^n$, 
    so the claim is true in this case. Otherwise, for $n = (\ell-1+j)$, $\phi_*\rho(x_0^n) = 
    (-1)^{j-1}(i_Z)_*s_j(N_{Z/X^*})_T =
    (-1)^{j-1}(i_Z)_*i_Z^*s_j(M)_T = 
    (-1)^{j-1}s_j(M)_Tc_{top}(M)_T$. Here $s_j(M)_T$ denotes the $j$-th 
    equivariant Segre class of the bundle $M = Q_{F_0}\oplus 
    \bigoplus_\lambda Q_{F_\lambda}\oplus\bigoplus 
    Q_{\Delta_i}\oplus\bigoplus_j Q_{\xi_j}$. The last equality follows from 
    lemma \ref{top Chern is Thom class}. Clearly this last expression is equal to 
    $\rho(P)$ for some universal polynomial $P$ in the Chern roots 
    of $M$. 
    We are therefore reduced to proving that $\nu_*(f) = P$. Let $s 
    =\mathrm{rank}(N_{W/Z})$ and let
    $\mu:\Bl{\C^{s+N}}\rightarrow \C^{s+N}$ denote the blow-up of 
    $\C^{s+N}$ along $\C^s\times 0$, where we give both spaces the structure 
    of toric varieties, with big torus $L$. $\C[\Bl{\Sigma}]$ 
    and $\C[\Sigma]$ both correspond to the rings of piece-wise 
    polynomial functions on the fans of $\Bl{\C^{s+N}}$ and 
    $\C^{s+N}$, which are in turn isomorphic to the $L$-equivariant 
    cohomology rings of the above spaces.
    We may view $\nu_*$ as the polyhedral version of the equivariant pushforward 
    $\mu_* :H^*_L(\Bl{\C^{s+N}})\rightarrow H^*_L(\C^{s+N})$. Then 
    $\nu_*f = P$ follows from the fact that 
    $\mu_*(c_1(E)_L)^{N-1+j} = P$, where here we evaluate $P$ at the 
    equivariant Chern classes of the coordinate hyperplanes. 
\end{proof}

\begin{rmk}\rm
As in section \ref{Poly}, the above theorems 
hold for piece-wise convergent power series.
\end{rmk}

\section{Blow-up Formula for Orbifold Elliptic Genus}\label{Blow Up Formula}

Let $X$ be a smooth projective variety with a holomorphic $G\times T$ action, $D = \sum_I \alpha_i D_i$ a $G\times T$-invariant $G$-normal crossing divisor with coefficients $\alpha_i < 1$, $g,h \in G$ a pair of commuting elements, and $X^{g,h}_\gamma$ a component of the common fixed point locus of $g$ and $h$. Recall the definition given in section \ref{Definitions} for the orbifold elliptic class of the pair $(X^{g,h}_\gamma,D)$. The goal in this section is to prove the following theorem:

\begin{thm}\label{orbifold chang of var}
    Let $f: \Bl{X}\rightarrow X$ be the blow-up along a smooth $G\times 
    T$-invariant subvariety $Z$ which has normal crossings with 
    respect to the components of $D$. Fix a commuting pair $g,h \in 
    G$ and a component $X^{g,h}_\gamma$ of the fixed point locus 
    $X^{g,h}$. Let $\set{\Bl{X}^{g,h}_\mu}$ denote the components of 
    $\Bl{X}^{g,h}$ which get mapped to $X^{g,h}_\gamma$ under $f$. 
    Then:
    \begin{align*}
	f_*\sum_\mu \mathcal{E}ll_{orb}(\Bl{X}^{g,h}_\mu,\Bl{D})
    &= \mathcal{E}ll_{orb}(X^{g,h}_\gamma,D)
    \end{align*}    
    Here $\Bl{D}$ is the divisor satisfying $f^*(K_X+D) = 
    K_{\Bl{X}}+\Bl{D}$.
\end{thm}

\begin{proof}
    By the projection formula, it suffices to assume that every 
    component of $D$ intersects $Z$ with multiplicity $1$. 
    Furthermore, by applying deformation to the normal cone (lemma 
    \ref{deformation}) we may assume that $X$ is the normal cone 
    space $X^*$ described in section \ref{Twisted Polyhedral Complex}. 
    Using the notation from section \ref{Twisted Polyhedral Complex}, 
    we have the following correspondences:
    \begin{align*}
	H &\longleftrightarrow (g,h)\\
	Z\cap X^{g,h}_\gamma &\longleftrightarrow W\\
	\set{D_j}_{j\in I^{g,h}_\gamma}&\longleftrightarrow \set{V_{\xi_j}}\\
	\set{D_i}_{i\not \in I^{g,h}_\gamma}&\longleftrightarrow 
	\set{V_{\Delta_i}}\\
	X^{g,h}_\gamma &\longleftrightarrow \bigcap_\ep V_{N_\ep}\cap\bigcap_{\lambda}
	V_{F_\lambda}\cap\bigcap_j V_{\xi_j}\\
	Z &\longleftrightarrow V_{F_0}\cap \bigcap_{\lambda}V_{F_\lambda}
	\cap \bigcap_{i}V_{\Delta_i}\cap \bigcap_j V_{\xi_j}
    \end{align*}
    
    Given these identifications, 
    $\mathcal{E}ll_{orb}(X^{g,h}_\gamma,D)=$
    \begin{align*}
	&\prod_{TX^*}\ellip{\frac{x_i}{2\pi i}}
	\prod_{\lambda,Q_{F_\lambda}}\frac{
	\theta(\twopi{x_{\lambda,c}})
	\theta(\twopi{x_{\lambda,c}}+\lambda(g)-\lambda(h)\tau-z)}
	{\theta(\twopi{x_{\lambda,c}}-z)
	\theta(\twopi{x_{\lambda,c}}+\lambda(g)-\lambda(h)\tau)}\times\\
	&\prod_{\ep,Q_{N_\ep}}\frac{
	\theta(\twopi{x_{\ep,a}})
	\theta(\twopi{x_{\ep,a}}+\ep(g)-\ep(h)\tau-z)}
	{\theta(\twopi{x_{\ep,a}}-z)
	\theta(\twopi{x_{\ep,a}}+\ep(g)-\ep(h)\tau)}\times\\
	&\prod_j\frac{\theta(\twopi{\xi_j})
	\theta(\twopi{\xi_j}+\lambda_j(g)-\lambda_j(h)\tau-(-\alpha_j+1)z)
	\theta(-z)}
	{\theta(\twopi{\xi_j}-z)
	\theta(\twopi{\xi_j}+\lambda_j(g)-\lambda_j(h)\tau)
	\theta(-(-\alpha_j+1)z)}\times\\
	&\prod_{i}\jacc{\frac{\Delta_i}{2\pi i}}{(-\alpha_i+1)}
    \end{align*}
    
    Here $x_{\lambda,c}, x_{\ep,a}, \xi_j$, and $\Delta_i$ denote the 
    equivariant Chern roots of the bundles,
    $Q_{F_\lambda}, Q_{N_\ep}, Q_{\xi_j}$, 
    and $Q_{\Delta_i}$, respectively. The above equality follows from 
    lemma \ref{top Chern is Thom class}.
    
    
    Let $f: \Bl{X^*}\rightarrow X^*$ be the blow-up along $Z$ with 
    exceptional divisor $E$. Before proceeding further, it will be 
    convenient to set up some new notation.
    Consider the collection $I = \set{N_\ep,F_0,F_\lambda,\xi_j,\Delta_i,E}$.
    Let $\Bl{Q}_{N_\ep} = f^*Q_{N_\ep}$. For $A \in I$, $A\neq N_\ep, 
    E$, let
    $\Bl{Q}_{A} = f^*Q_{A}\otimes\mathcal{O}(-E)$. For ease of 
    notation later on, we also define $\Bl{Q}_E = \mathcal{O}(E)$. 
    For $A \in I-\set{E}$, we let $\Bl{V}_A$ denote the proper 
    transform of $V_A$, and we let $\Bl{V}_E$ simply equal $E$.
     
    Let $\Bl{X^*}^{g,h}_\mu$ be 
    a connected component of $\Bl{X^*}^{g,h}$ which gets mapped to 
    ${X}^{g,h}_\gamma$ under $f$. This space is the complete 
    intersection of subvarieties $\Bl{V}_A$ for $A$ in some indexing 
    set $I_\mu \subset I$. For $A \in I_\mu$, 
    $\Bl{Q}_A|_{\Bl{X^*}^{g,h}_\mu} = N_A\otimes\Bl{\lambda}_A$ for some 
    $(g,h)$ character $\Bl{\lambda}_{A,\mu}$ and some $(g,h)$-trivial 
    bundle $N_{A,\mu}$. We extend the definition of 
    $\Bl{\lambda}_{A,\mu}$ by 
    letting $\Bl{\lambda}_{A,\mu} = 0$ for $A \not \in I_\mu$.
    
    Finally, we define indices $\beta_{A}$ as follows: For $A = 
    \xi_j,\Delta_i,$ or $E$, we let $\beta_A$ be the coefficient of 
    $\Bl{V}_A$ as a divisor in $\Bl{D}$. Otherwise, we let $\beta_A = 
    0$.
    
    Applying 
    corollary \ref{blow up Chern class}, lemma \ref{blow up Thom 
    class}, and the above definitions,
    we obtain the following convenient
    expression for $\mathcal{E}ll_{orb}(\Bl{X^*}^{g,h}_\mu,\Bl{D})$:
    
    \begin{align*}
	&f^*\Big\{\prod_{TX^*}\ellip{\twopi{x_i}}
	\prod_{m,Q_A,A\in I}\ellipinv{\twopi{a_m}}\Big \}\times\\
	&\prod_{\ell,\Bl{Q}_A,A\in I}\frac{\twopi{\tilde{a}_\ell}
	\theta(\twopi{\tilde{a}_\ell}+\Bl{\lambda}_{A,\mu}(g)-\Bl{\lambda}_{A,\mu}(h)\tau
	-(-\beta_A+1)z)\theta'(0)}
	{\theta(\twopi{\tilde{a}_\ell}+\Bl{\lambda}_{A,\mu}(g)-\Bl{\lambda}_{A,\mu}(h)\tau)
	\theta(-(-\beta_A+1)z)}
    \end{align*}
    
    For $A \in I-\set{E}$ we define, for convenience of notation, 
    $\lambda_A\in R((g,h))$ and $\alpha_A\in \Q$ as follows: For $A = 
    N_\ep, F_\lambda, \xi_j$, and $\Delta_i$, we let $\lambda_A = \ep, 
    \lambda, \lambda_j$, and $\lambda_i$, respectively. Otherwise, we 
    set $\lambda_A = 0$. Next, for $A = \xi_j$ or $\Delta_i$, we 
    let $\alpha_A = \alpha_j$ or $\alpha_i$. Otherwise, we set 
    $\alpha_A = 0$.
    
    By the projection formula, we are reduced to proving the following 
    formula:
    
    \begin{align}
	&f_*\sum_{\mu}\label{exp 1}
	\prod_{\ell,\Bl{Q}_A,A\in I}\frac{\twopi{\tilde{a}_\ell}
	\theta(\twopi{\tilde{a}_\ell}+\Bl{\lambda}_{A,\mu}(g)-\Bl{\lambda}_{A,\mu}(h)\tau
	-(-\beta_A+1)z)\theta'(0)}
	{\theta(\twopi{\tilde{a}_\ell}+\Bl{\lambda}_{A,\mu}(g)-\Bl{\lambda}_{A,\mu}(h)\tau)
	\theta(-(-\beta_A+1)z)} =\\
	&\prod_{m,Q_A,A\in I-E}\frac{\twopi{a_m}\label{exp 2}
	\theta(\twopi{a_m}+\lambda_{A}(g)-\lambda_{A}(h)\tau
	-(-\alpha_A+1)z)\theta'(0)}
	{\theta(\twopi{a_m}+\lambda_{A}(g)-\lambda_{A}(h)\tau)
	\theta(-(-\alpha_A+1)z)}
    \end{align}
    
    Naturally, $\tilde{a}_\ell=\tilde{a}_\ell(A)$ denote the equivariant Chern roots of 
    $\Bl{Q}_A$ and $a_m = a_m(A)$ denote the equivariant Chern roots 
    of $Q_A$.
    
    Referring to the notation from section \ref{Twisted Polyhedral 
    Complex} for $A = N_\ep, F_0,F_\lambda,\Delta_k,\xi_j$ and $i=1,\ldots,\mathrm{rank}(A)$, 
    let us define $x_{A,i} \in \C[\Sigma] = w_{\ep,i}, 
    x_i,y_{\lambda,i},d_k, e_j$, 
    respectively. For $A \in I$ and $i=1,\ldots,\mathrm{rank}(A)$ we 
    define $\tilde{x}_{A,i} \in \C[\Bl{\Sigma}]$ similarly. Define $F \in 
    \C[|\Sigma|]$ to be the power series:
    \begin{align*}
	F = \prod_{i,A\in I-E}\frac{\twopi{x_{A,i}}
	\theta(\twopi{x_{A,_i}}+\lambda_{A}(g)-\lambda_{A}(h)\tau
	-(-\alpha_A+1)z)\theta'(0)}
	{\theta(\twopi{x_{A,i}}+\lambda_{A}(g)-\lambda_{A}(h)\tau)
	\theta(-(-\alpha_A+1)z)}
    \end{align*}
    
    Define $F_\mu \in \C[|\Bl{\Sigma}|]=$
    \begin{align*}
    \prod_{i,A\in I}\frac{\twopi{\tilde{x}_{A,i}}
	\theta(\twopi{\tilde{x}_{A,i}}+\Bl{\lambda}_{A,\mu}(g)-\Bl{\lambda}_{A,\mu}(h)\tau
	-(-\beta_A+1)z)\theta'(0)}
	{\theta(\twopi{\tilde{x}_{A,i}}+\Bl{\lambda}_{A,\mu}(g)-\Bl{\lambda}_{A,\mu}(h)\tau)
	\theta(-(-\beta_A+1)z)}
    \end{align*}
    
    Clearly $F$ is $\rho$ applied to expression \ref{exp 1} and 
    $F_\mu$ is $\rho$ applied to the $\mu$-th summand in expression 
    \ref{exp 2}. By lemma \ref{pushforward commutes}, we are reduced 
    to proving $\nu_* \sum_\mu F_\mu = F$. To do this, think of 
    $(\Sigma,N_\Sigma)$ as the polyhedral complex associated to the toric 
    variety $\C^M$ where $M = \dim \Sigma$. Let $\Bl{\C^M}$ be the 
    toric blow-up of $\C^M$ which corresponds to the polyhedral 
    subdivision $\Bl{\Sigma}\rightarrow \Sigma$ described in section 
    \ref{Twisted Polyhedral Complex}. We may view $g$ and $h$ as 
    elements of the big torus of $\C^M$, i.e., as elements of a 
    finite index sup-lattice of $N_\Sigma$. Under this 
    identification, $x_{A,i}(g) = \lambda_A(g)$.
    
    The $(g,h)$-fixed components of $\Bl{\C^M}$ are in one-one 
    correspondence with the fixed components $\Bl{X^*}^{g,h}_\mu$ and 
    in one-one correspondence with subcones $C_\mu \subset 
    \Bl{\Sigma}$. These are the cones of maximal dimension which 
    correspond to affine open sets $U_{C_\mu}$ of the form 
    $\C^a\times(\C^*)^b$, where the characters of the 
    $(g,h)$-representation $\C^a$ are all non-trivial. For $C\supset 
    C_\mu$, let $I(C)$ index the collection of piece-wise linear 
    functions $\tilde{x}_{A,i}$ which are dual to $C$. Since $g$ and 
    $h \in N_\Sigma$ lie inside $C_\mu$, it makes sense to evaluate 
    $\tilde{x}_{A,i}|_C$ at $g$ and $h$. One sees easily in fact 
    that $\tilde{x}_{A,i}|_C(g)=\tilde{\lambda}_{A,\mu}(g)$ and 
    similarly for $\tilde{x}_{A,i}|_C(h)$. When distinguishing between 
    different cones, it will be convenient to denote the 
    collection $\set{\tilde{x}_{A,i}|_C}_{I(C)}$ by 
    $\set{\tilde{x}_{C,j}}_{j=1}^{|I(C)|}$. With this notation, if 
    $\tilde{x}_{C,j} = \tilde{x}_{A,i}|_C$, we also define 
    $\beta_j = \beta_A$. We define $x_{C,j}$ and $\alpha_j$ similarly 
    when $C\subset \Sigma$. 
      
    From this it follows that for $C\supset C_\mu$:
    
    \begin{align*}
	F_\mu|_{C}= 
	\prod_{j=1}^{\dim C}\frac{\twopi{\tilde{x}_{C,j}}
	\theta(\twopi{\tilde{x}_{C,j}}+\tilde{x}_{C,j}(g)-\tilde{x}_{C,j}(h)\tau
	-(-\beta_j+1)z)\theta'(0)}
	{\theta(\twopi{\tilde{x}_{C,j}}+\tilde{x}_{C,j}(g)-\tilde{x}_{C,j}(h)\tau)
	\theta(-(-\beta_j+1)z)}
    \end{align*}
    
    Otherwise, if $C_\mu$ is not contained in $C$, it is easy to see 
    that $F_\mu|_C = 0$.
    
    Now let $C_\gamma \subset \Sigma$ be the cone which corresponds 
    to $(X^*)^{g,h}_\gamma$. Fix a cone $K\subset \Sigma$ containing 
    $C_\gamma$. Let $\Bl{\Sigma}_K$ denote the subdivision of $K$ inside 
    $\Bl{\Sigma}$. Each cone $C \subset \Bl{\Sigma}_K$ with the 
    same dimension as $K$ contains a unique cone $C_\mu$ for some 
    $\mu$. Moreover, every cone $C_\mu$ is contained in one such cone 
    $C\subset \Bl{\Sigma}_K$. We therefore have that
    
    \begin{align*}
	\sum_\mu F_\mu|_{C} =
	\prod_{j=1}^{\dim K}\frac{\twopi{\tilde{x}_{C,j}}
	\theta(\twopi{\tilde{x}_{C,j}}+\tilde{x}_{C,j}(g)-\tilde{x}_{C,j}(h)\tau
	-(-\beta_j+1)z)\theta'(0)}
	{\theta(\twopi{\tilde{x}_{C,j}}+\tilde{x}_{C,j}(g)-\tilde{x}_{C,j}(h)\tau)
	\theta(-(-\beta_j+1)z)}
    \end{align*}
    
    To complete the proof, it remains to show that $(\nu_*\sum 
    F_\mu)|_K = F|_K$. By the push-forward formula for $\nu_*$, this 
    is equivalent to proving the identity:
    
    \begin{align*}
	&\sum_{C\subset \Bl{\Sigma}_K}
	\prod_{j=1}^{\dim K}\frac{
	\theta(\twopi{\tilde{x}_{C,j}}+\tilde{x}_{C,j}(g)-\tilde{x}_{C,j}(h)\tau
	-(-\beta_j+1)z)\theta'(0)}
	{\theta(\twopi{\tilde{x}_{C,j}}+\tilde{x}_{C,j}(g)-\tilde{x}_{C,j}(h)\tau)
	\theta(-(-\beta_j+1)z)}=\\
	&\prod_{j=1}^{\dim K}\frac{
	\theta(\twopi{x_{K,j}}+x_{K,j}(g)-x_{K,j}(h)\tau
	-(-\alpha_j+1)z)\theta'(0)}
	{\theta(\twopi{x_{K,j}}+x_{K,j}(g)-x_{K,j}(h)\tau)
	\theta(-(-\alpha_j+1)z)}
   \end{align*}
   
   Here the functions $\tilde{x}_{C,j}$ are regarded as linear 
   combinations of the functions $x_{K,j}$. The above formula follows 
   from theorem $7$ of the preprint \cite{RW} or by lemma $8.1$ 
   of \cite{BL}. This completes the proof.
\end{proof}


Now, let $Z$ be a projective $\Q$-Gorenstein variety with log terminal singularities, and a regular $G\times T$ action.  
Let $f: X \rightarrow Z$ and $g: Y\rightarrow Z$ be two equivariant resolutions of singularities. We assume that the exceptional locus of both resolutions is a $G$-normal divisor with simple normal crossings. Define $D_X$ so that $K_X + D_X = f^*K_Z$; define $D_Y$ similarly. Then the equivariant orbifold elliptic genera of $(X,D_X)$ and $(Y,D_Y)$ coincide. Indeed, by the equivariant version of the weak factorization theorem \cite{W}, we may connect $X$ to $Y$ by a sequence of equivariant blow-ups and blow-downs in such a way that the blow-ups at each intermediate pair $(X_i,D_{X_i})$ occur along a smooth base with normal crossings with respect to the components of $D_{X_i}$. Moreover, the procedure described in \cite{BL} theorem $3.7$ to make the intermediate pairs $(X_i,D_{X_i})$ $G$-normal extends to the $T$-equivariant case. Hence, by the equivariant change of variables formula, $Ell_{orb}(X,D_X,G) = Ell_{orb}(Y,D_Y,G)$.

We will, however, require a slightly stronger version of the above result for the purposes of this paper:

\begin{thm}
Let $f: (X,D_X) \rightarrow (Y,D_Y)$ be a $G\times T$-equivariant birational morphism between smooth, equivariant, $G$-normal log terminal pairs $(X,D_X)$ and $(Y,D_Y)$. Assume furthermore that $f^*(K_Y+D_Y) = K_X + D_X$. Then $f_*\mathcal{E}ll_{orb}(X,D_X,G) = \mathcal{E}ll_{orb}(Y,D_Y,G)$.
\end{thm}

\begin{proof}
The weak factorization theorem allows us to factor $f$ into a sequence of equivariant blow-ups and blow-downs

$$X=X_0\dashrightarrow X_1 \dashrightarrow\cdots\dashrightarrow X_k = Y$$
such that for some intermediate index $i_0$ the maps $X_i \rightarrow X$ are morphisms for $i \leq i_0$ and the maps $X_i \rightarrow Y$ are morphisms for $i \geq i_0$. Moreover, by \cite{BL} theorem $3.7$, we may still guarantee that all the intermediate varieties are $G$-normal with respect to the appropriate divisors. Note that since $f$ is itself a smooth morphism, we may conclude that the maps $X_i\rightarrow Y$ are smooth morphisms for all i.

We now apply induction on $k$. For $k = 1$, the theorem is obvious. Otherwise, consider the intermediate variety $X_1$. By assumption, $X_1 \neq X,Y$. Either $X \dashrightarrow X_1$ is a blow-up or blow-down. Suppose first that $X \leftarrow X_1$ is a blow-down. Call this morphism $g$. Define $D_1$ so that $K_{X_1}+D_1 = g^*(K_X+D_X)$. Note that if $h : X_1\rightarrow Y$ is the morphism $f\circ g$, then $K_{X_1}+D_1 = h^*(K_Y+D_Y)$. By the change of variables formula, $\mathcal{E}ll_{orb}(X,D_X,G) = g_*\mathcal{E}ll_{orb}(X_1,D_1,G)$. Therefore $f_*\mathcal{E}ll_{orb}(X,D_X,G) = h_*\mathcal{E}ll_{orb}(X_1,D_1,G)$. But this is equal to $\mathcal{E}ll_{orb}(Y,D_Y,G)$ by the induction hypothesis.

The case in which $X \rightarrow X_1$ is a blow-down is proved similarly. 
\end{proof}

\section{The Equivariant McKay Correspondence}\label{Equiv McKay}

Here we use the results from the preceding sections to prove an equivariant analogue of the McKay Correspondence for the elliptic genus. As a corollary, we will arrive at an equivariant version of the DMVV formula.

Let $X$ be a projective variety with a $G\times T$ action, where $G$ is a finite group and $T$ is a compact torus. Let $\psi:X\rightarrow X/G$ be the quotient morphism. Assume that $X/G$ has an equivariant crepant resolution $V$ and that $\psi^*(K_{X/G}) = K_X$. Let $F \subset X$ be a fixed component of the $G\times T$-action, and let $\set{P}\subset V$ denote the fixed components in $V$ which get mapped to $\psi(F)$. Then:

\begin{thm}\label{Local McKay}
$$\int_F \frac{\mathcal{E}ll_{orb}(X,G)}{e(F)}= \sum_P\int_P \frac{\mathcal{E}ll(V)}{e(P)}.$$
\end{thm}

\begin{proof}
Let $Z \rightarrow V$ be a sequence of equivariant blow-ups of $V$ so that the exceptional locus of the resolution $\pi:Z\rightarrow X/G$ is a divisor with simple normal crossings. Let $D_Z$ be the divisor on $Z$ such that $K_Z+D_Z = \pi^*K_{X/G}$. Define $\hat{Z}_0$ to be the normalization of $Z$ in the function field of $X$. By Abhyankar's Lemma, the induced map $\mu_0:\hat{Z}_0\rightarrow Z$ is a toroidal morphism of toroidal embeddings. Let $\hat{Z}$ be a projective toroidal resolution of singularities of $\hat{Z}_0$, and define $D_{\hat{Z}}$ so that $K_{\hat{Z}}+D_{\hat{Z}} = \mu^*(K_Z+D_Z)$, where $\mu:\hat{Z}\rightarrow Z$ is the obvious map. We may further assume that $\hat{Z}$ has a $G$ action, and that the pair $(\hat{Z},D_{\hat{Z}})$ is $G$-normal. (see \cite{AW}). We obtain the following commutative diagram:

$$\begin{CD}
\hat{Z} @>\mu>> Z \\
@V\phi VV 		  @VV\pi V \\
X @>\psi >> X/G \\
\end{CD}$$

Here, the vertical arrows are resolutions of singularities, and the horizontal arrows are birational to a quotient by $G$. It is clear that all the constructions involved (normalization, blow-up along a $T$-invariant ideal sheaf) are $T$-equivariant, and that consequently the above morphisms are $T$-equivariant.

Since $\phi^*(K_X) = K_{\hat{Z}}+D_{\hat{Z}}$, we have that $\mathcal{E}ll_{orb}^T(X,G) = \phi_*\mathcal{E}ll_{orb}({\hat{Z}},D_{\hat{Z}},G)$. Let $\set{L} \subset \hat{Z}$ denote the fixed components of $\hat{Z}$ which map to $F$. By functorial localization:

$$\int_F \frac{\mathcal{E}ll_{orb}(X,G)}{e(F)}= \sum_L\int_L \frac{\mathcal{E}ll_{orb}(\hat{Z},D_{\hat{Z}},G)}{e(L)}.$$

Now, by theorem \ref{Toroidal McKay}, $\mu_*\mathcal{E}ll_{orb}(\hat{Z},D_{\hat{Z}},G) = \mathcal{E}ll(Z,D_Z)$. Let $\set{K}\subset Z$ denote the fixed components which get mapped to $\psi(F)$ under the resolution $Z\rightarrow X/G$. Clearly $\set{L} = \phi^{-1}F = \phi^{-1}\psi^{-1}\psi(F) = \mu^{-1}\pi^{-1}\psi(F) = \mu^{-1}\set{K}$. Thus, by functorial localization applied to $\mu_*$, we have:

$$\sum_L\int_L\frac{\mathcal{E}ll_{orb}(\hat{Z},D_{\hat{Z}},G)}{e(L)}
= \sum_K\int_K\frac{\mathcal{E}ll(Z,D_Z)}{e(K)}.$$

Finally, since $\set{K}$ denotes the fixed components of $Z$ which get mapped to $\set{P}$, functorial localization applied to $Z \rightarrow V$ completes the proof.
\end{proof}

We now discuss some corollaries of the above result. We begin with a proof of the equivariant DMVV conjecture. Let $\Proj_2^{(n)} = (\Proj_2)^n/S_n$ denote the $n$th symmetric product of the projective plane. The natural group action of $T = S^1\times S^1$ on $\Proj_2$ extends to $(\Proj_2)^n$ in the obvious manner, and commutes with the action of $S_n$. The action of $T$ on $\Proj_2$ also extends to $\Proj_2^{[n]}$, the Hilbert scheme of $n$ points on $\Proj_2$, and the Hilbert-Chow morphism $\Proj_2^{[n]}\rightarrow \Proj_2^{(n)}$ is an equivariant crepant resolution.

Sitting inside $\Proj_2^{(n)}$ is the open variety $(\C^2)^{(n)}$. It has a single $T$-fixed point $p$, which is the image under the quotient morphism of the $S_n\times T$ fixed point $(0,0)\times...\times (0,0) \in (\C^2)^n$. The pre-image of $(\C^2)^{(n)}$ under the Hilbert-Chow morphism is just $(\C^2)^{[n]}$. Hence, the pre-image of $p$ under the Hilbert-Chow morphism is simply the collection of $T$-fixed points of $(\C^2)^{[n]}$. If $(u_1,u_2)$ denote the equivariant parameters of the $T$-action, let $t_j = e^{2\pi iu_j}$. Then the above theorem implies that:

$$Ell_{orb}((\C^2)^n,S_n,t_1,t_2) = Ell((\C^2)^{[n]},t_1,t_2).$$
Note that the LHS involves equivariant data associated to the single fixed point $p$, whereas the RHS involves a sum of equivariant data associated to partitions of $n$.

For $z$ the complex parameter appearing in the definition of the elliptic class and $\tau$ the lattice parameter used in the definition of the Jacobi theta function, let $y = e^{2\pi iz}$ and $q = e^{2\pi i\tau}$. Viewing the equivariant elliptic indices as formal power series in the variables $q,y,t_1,$ and $t_2$, and applying theorem $3.1$ of \cite{LLJ}, we have the following equivariant analogue of the DMVV formula:

\begin{thm}
Write $Ell(\C^2,t_1,t_2) = \sum_{m\geq 0,\ell,k}c(m,\ell,k)q^m y^\ell t_1^{k_1}t_2^{k_2}$. Then

$$\sum_{n>0} p^n Ell((\C^2)^{[n]},t_1,t_2) = 
\prod_{m\geq 0,n>0,\ell,k}\frac{1}{(1-p^n q^m y^\ell t_1^{k_1} t_2^{k_2})^{c(nm,\ell,k)}}$$
\end{thm}

We next discuss the equivariant elliptic genus analogue of the classical McKay correspondence, which was originally proved in \cite{RW}. Let $G \subset SU(2)$ be a finite subgroup. $G$ acts on $\C^2$ in the obvious fashion, and commutes with the diagonal action of $T = S^1$. Let $V \rightarrow \C^2/G$ be the crepant resolution of singularities. The action of $T$ lifts to $V$, and the fixed components of this action are compact. 

\begin{thm}
$Ell_{orb}(\C^2,G,t) = Ell(V,t)$
\end{thm}

\begin{proof}
View $\C^2$ as an open subset of $\Proj^2$, with the action of $G$ and $T$ extending to $\Proj^2$ in the obvious manner. The space $\Proj^2/G$ still has only an isolated singularity at the image of the origin $[0:0:1]$. Hence $\Proj^2/G$ has an equivariant crepant resolution which is a compactification of $V$. The above theorem then follows by letting $F = (0,0)$ in theorem \ref{Local McKay}.
\end{proof}

\section{Appendix}
\begin{lem}
   Let $f:X \rightarrow Y$ be a $T$-map of smooth compact simply 
    connected Kahler varieties. 
    Let $D \subset Y$ be a $T$-invariant divisor and let $E_i$ be 
    $T$-invariant normal crossing divisors on $X$ such that $f^*D = 
    \sum a_i E_i$ as Cartier divisors. Then for any 
    $\varepsilon$-regular neighborhood $U_{\varepsilon}$ of $D$ 
    there exist generators $\Theta_{E_i}^T$ for $c_1^T(E_i)$ and 
    $\Theta_D^T$ for $c_1^T(D)$ with the following properties:
    
    $(1)$ $\Theta_D^T$ has compact support in $U_{\varepsilon}$ and 
    $\Theta_{E_i}^T$ have compact support in $f^{-1}(U_{\varepsilon})$.
    
    $(2)$ $f^* \Theta_D^T = \sum a_i \Theta_{E_i}^T + d_T (\eta)$ on 
    the level of forms, where $\eta$ is a $T$-invariant form with 
    compact support in $U_{\varepsilon}$.
    
    $(3)$ $\Theta_D^T$ and $\Theta_{E_i}^T$ represent the equivariant 
    Thom classes of the varieties $D$ and $E_i$
\end{lem}

The only real issue above is to ensure that $\eta$ has compact support 
in the desired neighborhood.

\begin{proof}
    We first solve this problem in the non-equivariant category. For 
    $V$ any Cartier divisor, denote by $L_V$ the line bundle it 
    induces. Let $U_{\varepsilon}$
    be a $T$-invariant tubular neighborhood of $D$ of radius 
    $\varepsilon$. Outside $U_{\frac{\varepsilon}{2}}$, the constant 
    function $1$ is a section of $L_D$. Define a metric $h_{far}$ in this 
    region by $h_{far} = \norm{1}^2 \equiv 1$. Let $h_{near}$ be a metric 
    inside $U_{\varepsilon}$. Piece the two metrics into a global 
    metric $h$ on $L_D$ using a partition of unity. The first Chern class 
    of $L_D$ is then represented by the form 
    $\Theta_D =\frac{i}{2\pi}\delbar\del \log h$. This form clearly has compact 
    support in $U_{\varepsilon}$.
    
    Let $U_{\varepsilon_i}$ be tubular neighborhoods of $E_i$. 
    Choose $\varepsilon_i$ small enough so that each of these 
    neighborhoods is contained in $f^{-1}U_{\varepsilon}$. Define 
    metrics $h_i$ on $E_i$ in a manner analogous to the above 
    construction of $h$. Clearly the forms $\Theta_{E_i} = 
    \frac{i}{2\pi}\delbar\del \log h_i$ have compact support in 
    $U_{\varepsilon_i}$ and represent the first Chern classes of $L_{E_i}$. 
    
    We have two natural choices for a metric on $f^*L_D$, namely 
    $f^*h$ and $h_1^{a_1}\cdots h_k^{a_k}$. Choose a smooth nonzero 
    function $\varphi$ so that $f^*h = \varphi h_1^{a_1}\cdots h_k^{a_k}$.
    Notice that $\varphi \equiv 1$ outside $f^{-1}U_{\varepsilon}$. 
    We have:
    
    $$\delbar\del\log f^*h = \delbar\del\log\varphi +
    \sum_i a_i\delbar\del\log h_i.$$
    
    But this implies that $f^*\Theta_D = \sum_i a_i\Theta_{E_i} 
    +\frac{i}{2\pi}\delbar\del\log\varphi$. If we let $d^c = 
    \frac{i}{4\pi}(\delbar-\del)$, we may write this last equation as:
    
    $$f^*\Theta_D = \sum_i a_i\Theta_{E_i}-dd^c\log\varphi.$$
    
    The form $\eta = -d^c\log\varphi$ clearly has compact support in 
    $f^{-1}U_{\varepsilon}$. It remains to argue that $\Theta_D$ and 
    $\Theta_{E_i}$ represent the Thom classes of $D$ and $E_i$. 
    Suppose $V$ is any effective Cartier divisor on a smooth simply connected 
    compact Kahler variety. Then $L_V$ is nontrivial as a holomorphic line 
    bundle and therefore $c_1(L_V) \neq 0$. This follows from the 
    fact that $H^1(X,\mathcal{O}) = 0$ and therefore that the map 
    $c_1: H^1(X,\mathcal{O}^*)\rightarrow H^2(X,\Z)$ is injective. 
    Furthermore, $c_1(L_V)$ is clearly a non-torsion class in 
    $H^2(X,\Z)$. Hence we can find an $\omega$ such that $\int_X 
    c_1(L_V)\wedge \omega \neq 0$. Since $c_1(L_V)$ is Poincar\'e 
    dual to $V$, we have $\int_X c_1(L_V)\wedge \omega = \int_V 
    \omega = \int_X \Phi_V \wedge \omega$ where $\Phi_V$ is the Thom 
    class of the normal bundle $N_V$ to $V$. However, if $\Theta_V$ is a 
    representative of $c_1(L_V)$ with compact support in $N_V$, 
    we must have $\Theta_V = a\Phi_V + d\psi$ for some 
    form $\psi$ with compact support in $N_V$. Therefore 
    $\int_X \Theta_V\wedge \omega = a\int_X \Phi_V\wedge \omega = 
    \int_X \Phi_V \wedge \omega \neq 0$. It follows that $a=1$. 
    Thus, $\Theta_D$ and $\Theta_{E_i}$ indeed represent the Thom 
    classes of $D$ and $E_i$. I should remark that it is well-known 
    that $\Theta_V$ is Poincar\'e dual to $V$. However, under weaker 
    conditions, a form that is Poincar\'e dual to a submanifold $V$ 
    may not necessarily
    represent the Thom class of $N_V$. For example, if $V$ where 
    homologous to zero, then $0$ would certainly be Poincar\'e dual 
    to $V$ but would not represent its Thom class, which would be 
    non-trivial.
    
    This completes the non-equivariant 
    portion of the proof.
    
    By averaging over the group $T$, we may assume that all the forms 
    above are $T$-invariant. For notational simplicity, let us assume 
    that $T = S^1$. Let $V$ be the vectorfield on $X$ induced by the 
    $T$-action. Let $g_i$ be the functions compactly supported in 
    $f^{-1}U_{\varepsilon}$ which satisfy the moment map equation 
    $i_V\Theta_{E_i} = dg_i$. Similarly, let $W$ be the vectorfield 
    on $Y$ defined by the $T$-action and define $g$ so that it 
    satisfies $i_V\Theta_D = dg$ and has support inside 
    $U_{\varepsilon}$. Note that since $f$ is $T$-equivariant, 
    $i_Vf^*\Theta_D = f^*i_W\Theta_D = f^*dg$. We 
    then have $d(g\circ f) = \sum_i a_i dg_i + i_Vd\eta = \sum_i a_i dg_i 
    -di_V\eta$. Hence $g\circ f = \sum_i a_i g_i -i_V\eta$. But this implies 
    that:
    
    $$f^*(\Theta_D+g) = \sum_i a_i(\Theta_{E_i}+g_i)+(d-i_V)\eta.$$
    
    But this is precisely the relation we wish to in the equivariant 
    cohomology.
\end{proof}

\begin{lem} \bf{(Excess Intersection Formula)} \it
    Let $X$ be a smooth compact variety with irreducible normal crossing 
    divisors $D_1,\ldots, D_k$. For $I \subset \set{1,\ldots,k}$ 
    denote by $X_{I,j}$ the $j$th connected component of 
    $\cap_{I}D_i$ and by $\Phi_{I,j}$ its Thom class. Fix irreducible 
    subvarieties $X_{I_1,j_1}$ and $X_{I_2,j_2}$. For $I_0 = I_1 \cup 
    I_2$, let $X_{I_0,j}$ be the irreducible components of 
    $X_{I_1,j_1} \cap X_{I_2,j_2}$. Then:
    
    $$\Phi_{I_1,j_1}\wedge \Phi_{I_2,j_2} = 
    \sum_{I_0,j}\Phi_{I_0,j}\prod_{I_1\cap I_2}\Phi_{i}.$$
\end{lem}

\begin{proof}
    Let $N_{I,j}$ be tubular neighborhoods of $X_{I,j}$ which are 
    disjoint for each indexing set $I$ and which satisfy 
    $N_{I,j} \subset N_{I',j'}$ for $X_{I,j}\subset X_{I',j'}$. If we choose 
    $\Phi_i$ to have compact support in a sufficiently small tubular 
    neighborhood of $D_i$, then $\prod_I \Phi_i$ will have compact 
    support in $\coprod_{j}N_{I,j}$. Moreover, the extension by zero 
    of $(\prod_I \Phi_i)|_{N_{I,j}}$ will represent the Thom class of 
    $X_{I,j}$ (see [BT]). We may also ensure that 
    $\Phi_{I_1,j_1}\wedge \Phi_{I_2,j_2}$ has compact support in 
    $\coprod_j N_{I_0,j}$. Thus:
    
    $$\Phi_{I_1,j_1}\wedge \Phi_{I_2,j_2} = \sum_{I_0,j}
    \big( \prod_{I_1}\Phi_i \prod_{I_2}\Phi_i \big)|_{N_{I_0,j}} =
    \sum_{I_0,j}\big(\prod_{I_0}\Phi_i \prod_{I_1\cap 
    I_2}\Phi_i\big)|_{N_{I_0,j}} = $$
    $$\sum_{I_0,j}\big(\prod_{I_0}\Phi_i \big)|_{N_{I_0,j}} \prod_{I_1\cap 
    I_2}\Phi_i.$$
    
    This yields the desired formula.
\end{proof}

\begin{rmk}
    \rm Note that the above proof clearly extends to the equivariant 
    category.
\end{rmk}


\end{document}